\documentclass[12pt]{article}
\usepackage[russian]{babel}
\usepackage{graphicx}
\usepackage[bf]{caption2}
\usepackage{amsfonts, amssymb, amsthm, amsmath}
\usepackage{color}
\newcounter{theorem}

\newcounter{theoremcounter}

\newcounter{corollarycounter}
\newtheorem{theorem}[theoremcounter]{Theorem}

\newtheorem{corollary}[corollarycounter]{Corollary}
\newcommand{\la}{\lambda}
\newcommand{\vp}{{\mathbf p}}
\newcommand{\vz}{{\mathbf z}}
\newcommand{\vf}{{\mathbf f}}
\newcommand{\vx}{{\mathbf x}}

\newcommand{\vy}{{\mathbf y}}
\newcommand{\vu}{{\mathbf u}}

\usepackage{enumerate}

\begin{document}

\title{On the Three Methods for Bounding the Rate of Convergence for some Continuous-time Markov Chains}

\author{Alexander Zeifman\thanks{Vologda State
University, Vologda, Russia; Institute of Informatics Problems,
Federal Research Center ``Computer Science and Control'',  Russian
Academy of Sciences, Moscow, Russia; Vologda Research Center of the
Russian Academy of Sciences, Vologda, Russia. E-mail:
a$\_$zeifman@mail.ru} \and  Yacov Satin\thanks{Vologda State
University, Vologda, Russia}\and  Anastasia Kryukova\thanks{Vologda
State University, Vologda, Russia}\and Rostislav
Razumchik\thanks{Institute of Informatics Problems Federal Research
Center ``Computer Science and Control'',  Russian Academy of
Sciences, Moscow, Russia; Peoples' Friendship University of Russia
(RUDN University), Moscow, Russia}\and  Ksenia
Kiseleva\thanks{Vologda State University, Vologda, Russia} \and
Galina Shilova\thanks{Vologda State University, Vologda, Russia}}

\date{}

\maketitle

{\bf Abstract.} Consideration is given to the three different
analytical methods for the computation of upper bounds for the rate
of convergence to the limiting regime of one specific class of
(in)homogeneous continuous-time Markov chains. This class is
particularly suited to describe evolutions of the total number of
customers in (in)homogeneous $M/M/S$ queueing systems with possibly
state-dependent arrival and service intensities, batch arrivals and
services. One of the methods is based on the logarithmic norm of a
linear operator function; the other two rely on Lyapunov functions
and differential inequalities, respectively. Less restrictive
conditions (compared to those known from the literature) under which
the methods are applicable, are being formulated. Two numerical
examples are given. It is also shown that for homogeneous
birth-death Markov processes defined on a finite state space with
all transition rates being positive, all methods yield the same
sharp upper bound.

\bigskip


\section{Introduction}
\label{intro}

In this paper we revisit the problem of finding the upper bounds
for the rate of convergence of (in)homogeneous continuous-time Markov chains.
Consideration is given to classic inhomogeneous birth-death processes
and to special inhomogeneous chains with transitions intensities,
which do not depend on the current state.
Specifically, let $\{X(t), \ t\geq 0 \}$ be an inhomogeneous continuous-time Markov chain
with the state space ${\mathcal{X}=\{ 0, 1, 2, \dots, S\}}$,
where $1\le S \le \infty$.
Denote by $p_{ij}(s,t)=P\left\{ X(t)=j\left| X(s)=i\right. \right\}$,
$i,j \ge 0, \;0\leq s\leq t$, the transition probabilities of
$X(t)$ and by  $p_i(t)=P \left\{ X(t) =i \right\}$ -- the
probability that $X(t)$ is in state $i$ at time $t$.
Let $\vp(t) = \left(p_0(t), p_1(t), \dots, p_S(t)\right)^T$ be
probability distribution vector at instant $t$.
Throughout the paper it is assumed that in a small time interval
$h$ the possible transitions and their associated probabilities are
\begin{equation*}
p_{ij}(t,t+h)=
\left\{ \!
\begin{array}{cc}
q_{ij}(t)  h \!+\! \alpha_{ij}\left(t, h\right), & \mbox {if } j\neq i,
\\ 1 \! - \! \sum\limits_{k \in \mathcal{X}, k\neq i}q_{ik}(t) h \!+\! \alpha_{i}\left(
t,h\right), & \mbox {if } j=i,
\end{array}
\right.
%  \label{4001}
\end{equation*}
\noindent
where transition intensities $q_{ij}(t) \ge 0$ are arbitrary\footnote{It
is not required (as, for example, in \cite{Zeifman2018c}),
that $q_{i,i+k}(t)$ and $q_{i,i-k}(t)$ are
monotonically decreasing in $k$ for any $t \ge 0$.} non-random functions of $t$,
locally integrable on $[0,\infty)$, satisfying $\sup_{i \in \mathcal{X}}
\left (\sum_{k \in \mathcal{X}, k\neq i}q_{ik}(t) \right )\le L < \infty$ for almost all $t \ge 0$,
and $|\alpha_i(t,h)| = o(h)$ for $S<\infty$
and $\sup_{i \in \mathcal{X}}  |\alpha_i(t,h)| = o(h)$ for $S=\infty$.
The results of this paper are applicable to Markov chains $X(t)$
with the following transition intensities:
\begin{enumerate}[i.]
  \item $q_{ij}(t)=0$ for any $t\ge 0$ if $|i-j|>1$ and both $q_{i,i+1}(t)$ and $q_{i,i-1}(t)$ may depend on $i$;
  \item $q_{i,i-k}(t)=0$ for $k \ge 2$, $q_{i,i-1}(t)$ may depend on $i$
  and $q_{i,i+k}(t)$, $k \ge 1$, depend only on $k$;
\item $q_{i,i-k}(t)=0$ for $k \ge 1$ depend only on $k$, $q_{i,i+1}(t)$ may depend on $i$
  and $q_{i,i+k}(t)=0$, $k \ge 2$;
  \item both $q_{i,i-k}(t)$ and $q_{i,i+k}(t)$, $k \ge 1$, depend only on $k$ and do not depend on $i$.
\end{enumerate}

\noindent Being motivated by the application of the obtained results
in the theory of queues\footnote{Yet the scope of the obtained
results is not limited to queueing systems and includes a number of
other stochastic systems which appear, for example, in medicine and
biology, which satisfy the adopted assumptions.}, in what follows it
is convenient to think of $X(t)$ as of the process describing the
evolution of the total number of customers of a queueing system.
Then type (i) transitions describe Markovian queues with possibly
state-dependent arrival and service intensities (for example, the
classic $M_t(n)/M_t(n)/1$ queue); type (ii) transitions allow
consideration of Markovian queues with state-independent batch
arrivals and state-dependent service intensity; type (iii)
transitions lead to Markovian queues with possible state-dependent
arrival intensity and state-independent batch service; type (iv)
transitions describe Markovian queues with state-independent batch
arrivals and batch service. For the details concerning possible
applications of Markovian queues with time-dependent transitions we
can refer to \cite{Schwarz2016}, which contains a broad overview and
a classification of time-dependent queueing systems considered up to
2016 and  also
\cite{DiCrescenzo2018,Giorno2014,Granovsky2004,Schwarz2016,Zeifman2006,Zeifman2014b,Vvedenskaya2018,Olwal2012,radoslav,LI2007,Almasi2005,Moiseev2016,Brugno2017,Trejo2019}
and references therein.

In this paper we propose three different analytical methods
for the computation of the upper bounds\footnote{I.e. bounds which guarantee
that after a certain time, say $t^*$, the probability characteristics of the process $X(t)$ do
not depend on the initial conditions (up to a given discrepancy). Since the proposed
methods are analytic we do not compare them here from the numerical point of view (i.e.
memory requirement, speed, running time etc.).} for the rate of convergence to the limiting regime
(provided that it exists) of any process $X(t)$ belonging to one of the classes (i)--(iv).
The {first one is based on {\it logarithmic norm} of a linear operator function.
The second one uses simplest {\it Lyapunov functions}
and the third one relies on the {\it differential inequalities}.
Even though the methods are not new, it is the first time when it is shown
how they can be applied for the analysis of Markov chains with the transition intensities
specified by (i)--(iv). This constitutes the main contribution of the paper.
Another contribution is the fact that in a case of periodic intensities
the bounds on the rate of convergence depend on the intensities only through their mean
values over one period.

It is worth noticing here that, except for  the upper bounds for the
rate of convergence, we may also be interested in the lower bounds,
stability (perturbation) bounds or truncation bounds (with error
estimation). But the exact estimates of the rate of convergence
yield exact estimates of stability bounds (see, for example,
\cite{Kartashov1985,Liu2012,Mitrophanov2003,Mitrophanov2004,Rudolf,Zeifman1985,Mitrophanov2018}
and references therein). Moreover, as our research shows (see
\cite{Zeifman2006,Zeifman2014a,Zeifman2014b,Zeifman2018c}), in some
cases, all these quantities can be constructed automatically, given
that some good upper bounds for the rate of convergence are
provided. This makes us believe that the upper bounds are of primary
interest.

Estimation of the convergence rate
by virtue of the methods proposed in this paper
heavily relies on the notion of the reduced intensity matrix, say $B(t)$,
of a Markov chain $X(t)$. The matrix $B(t)$ can be obtained by
considering the probabilistic dynamics of the process $X(t)$,
given by the forward Kolmogorov system
\begin{equation}
\label{ur01}
\frac{d}{dt}\vp(t)=A(t)\vp(t),
\end{equation}
\noindent where $A(t)$ is the transposed intensity matrix i.e.
$a_{ij}(t)=q_{ji}(t)$, $i,j \in \mathcal{X}$. Since due to the
normalization condition $p_0(t) = 1 - \sum_{i = 1}^S p_i(t)$, we
can rewrite\footnote{For the detailed discussion of the
transformation \eqref{2.06} see \cite{Granovsky2004,Zeifman2006}.}
the system (\ref{ur01}) as follows:
\begin{equation}
\frac{d}{dt}{\vz}(t)= B(t){\vz}(t)+{\vf}(t), \label{2.06}
\end{equation}
\noindent where
$$
{\vf}(t)=\left( a_{10}(t),  a_{20}(t),\dots \right)^{T}, \
{\vz}(t)=\left(p_1(t), p_2(t),\dots \right)^{T},
$$
\begin{equation}
{\footnotesize
\! B(t)\!=\!
%\!=\! \left(b_{ij}(t)\right)_{i,j=1}^{\infty} \!=\!
\left(
\begin{array}{ccccc}
a_{11} \!-\!a_{10} & a_{12}
\!-\!a_{10} & \cdots & a_{1r}
\!-\!a_{10} & \cdots \\
a_{21} \!-\!a_{20} & a_{22}
\!-\!a_{20} & \cdots & a_{2r}
\!-\!a_{20} & \cdots \\
\cdots & \cdots & \cdots & \cdots  & \cdots \\
a_{r1} \!-\!a_{r0} & a_{r2}
\!-\!a_{r0} & \cdots & a_{rr}
\!-\!a_{r0} & \cdots \\
\vdots & \vdots & \vdots & \vdots  & \ddots
\end{array}
\right)\!.}\label{2.07}
\end{equation}

\noindent Here and henceforth each entry of $B(t)$ may depend on $t$
but, for the sake of brevity, the argument is omitted.
We note that the matrix $B(t)$  has no any probabilistic meaning.
All bounds of the rate of convergence to the limiting regime for
$X(t)$ correspond to the same bounds of the solutions of the system
\begin{equation}
\frac{d}{dt}{\vy}(t)= B(t){ \vy}(t), \label{hom1}
\end{equation}
\noindent because ${ \vy}(t) =
{\vz}^{*}(t)-{\vz}^{**}(t)$ is the difference of two solutions of
system (\ref{2.06}), and  ${ \vy}(t) = \left({y}_1(t), {y}_2(t),
\dots, {y}_S(t)\right)^T$ is the vector with the coordinates of
arbitrary signs. As it was firstly noticed in \cite{Zeifman1989}, it
is more convenient to study the rate of convergence using the
transformed version $B^*(t)$ of $B(t)$ given by
$B^*(t)=TB(t)T^{-1}$, where $T$ is the $S \times S$ upper triangular
matrix of the form
\begin{equation}
T=\left(
\begin{array}{ccccccc}
1   & 1 & 1 & \cdots & 1 \\
0   & 1  & 1  &   \cdots & 1 \\
0   & 0  & 1  &   \cdots & 1 \\
\vdots & \vdots & \vdots & \ddots \\
0   & 0  & 0  &   \cdots & 1
\end{array}
\right). \label{vspmatr}
\end{equation}

\noindent Let ${\vu}(t)=T{\vy}(t)$. Then the system \eqref{hom1} can be rewritten in the form
\begin{equation}
\frac{d}{dt} {\vu}(t) = B^*(t) {\vu}(t), \label{hom11}
\end{equation}

\noindent where ${\vu}(t) = \left( {u}_1(t), u_2(t), \dots, u_S(t)\right)^T$
is the vector with the coordinates of arbitrary signs. If one of the two matrices  $B^*(t)$ or $B(t)$
is known, the other is also (uniquely) defined.

The method based on the logarithmic norm of a linear operator
function and the corresponding bounds for the Cauchy operator of
reduced forward Kolmogorov system has already been applied
successfully in many settings (see, for example,
\cite{Granovsky2004,Zeifman2006}). Moreover, in \cite{Zeifman2018c}
we have obtained the bounds for the rate of convergence and
perturbation bounds for a process $X(t)$ belonging to classes
(i)--(iv) under the assumption that $B^*(t)$ is essentially
non-negative i.e. $b^*_{ij}(t) \ge 0$, $i \neq j$, $i,j \in
\mathcal{X}\backslash (0)$. The obtained bounds are sharp for
non-negative difference of the initial probability distributions of
$X(t)$. In this paper it is no longer assumed that $B^*(t)$ must be
essentially non-negative. Thus the considered class of eligible
processes $X(t)$ is wider than the one considered in
\cite{Zeifman2018c}.

It may happen that the difference of the initial probability
distributions of $X(t)$ has coordinates of different signs and/or
$B^*(t)$ contains negative elements. In such situations the upper
bounds provided by the method based on the logarithmic norm may not
be sharp. Having alternative estimates, provided by the other two
methods considered in this paper, we can choose the best. The idea
to apply Lyapunov functions for the analysis of Markov chains is not
new\footnote{For the detailed description of the approach we can
also refer to \cite{Meyn1993,Meyn2012}.} (see, for example,
\cite{Kalash1971,Malysh1982}). Yet, to our best knowledge, in the
considered setting they have not been applied before (see the
seminal paper \cite{Zeifman2019icumt}). The approach based on the
differential inequalities (see \cite{Zeifman2019rodos}) seems to be
the most general: it can be applied both in the case when $B(t)$ is
essentially non-negative (and yield the same results as the method
based on the logarithmic norm) and in other cases in which the other
two methods are not applicable. Usually the three methods lead to
different upper bounds and the quality (sharpness) of the bounds
depends on the properties of $B^*(t)$. All three methods are
applicable when the state space i.e. $S<\infty$. For countable
$\mathcal{X}$ the method based on Lyapunov functions no longer
applies. Note also that for a $X(t)$ with a finite state space
belonging to classes (i)--(iv) apparently no general method for the
construction of Lyapunov functions can be suggested. Thus here
consideration is given only to such $X(t)$, for which it can be
guessed how Lyapunov functions can be constructed.

The paper is structured as follows.
In the next section the explicit forms of the reduced intensity matrix
$B^*(t)$ for each class (i)--(iv) are given.
In section~3 we review the upper bounds on the rate of convergence,
obtained by the method based on the logarithmic norm.
Alternative upper bounds provided by Lyapunov functions and differential inequalities
for some $X(t)$ from classes (i)--(iv) are given in sections~4 and 5. Section~6 concludes the paper.

%the general description
%of the system under consideration is given
%and all the necessary notation is introduced.
%Section 3 contains the main result of the paper
%i.e. the theorem which specifies the
%convergence bounds.
%Section 4 provides explicit expressions
%for functions needed to compute
%convergence bounds.
%In the last two sections we provide
%extensive numerical examples and
%give some directions of
%further research.

\section{Explicit forms of the reduced intensity matrix}

As it was mentioned above,
estimation of the convergence rate of $X(t)$ to the limiting regime
is based on the reduced intensity matrix $B(t)$, given by \eqref{2.07}, or its transform $B^*(t)=TB(t)T^{-1}$.
In this section the explicit form of $B^*(t)$ for each class (i)--(iv) is given.

\subsection{$B^*(t)$ for $X(t)$ belonging to class $\rm (i)$}

Consider a process $X(t)$ with $a_{ij}(t)=0$ for any $t\ge 0$ if $|i-j|>1$,
$a_{i,i+1}(t)=\mu_{i+1}(t)$ and $a_{i+1,i}(t)=\la_i(t)$.
Then $X(t)$ is the inhomogeneous birth-death process with state-dependent
transition intensities $\la_i(t)$ (birth) and $\mu_{i+1}(t)$ (death).
In the queueing theory context, $X(t)$ describes the evolution
of the total number of customers in the
$M_n(t)/M_n(t)/1/S$ queue. For such
$X(t)$ in the case of countable state space (i.e. $S=\infty$)
the matrix $B^*(t)$ has the form:
$$B^*(t) =  $$
\begin{equation}
{\tiny \left(
\begin{array}{ccccccc}
-\left(\lambda_0+\mu_1\right)  & \mu_1
 & 0 & \cdots & 0 &\cdots &\cdots\\
\lambda_1  & -\left(\lambda_1+\mu_2\right) & \mu_2 & \cdots & 0 &\cdots &\cdots\\
\ddots & \ddots & \ddots & \ddots & \ddots  &\cdots\\
0 & \cdots & \cdots & \lambda_{r-1} &
-\left(\lambda_{r-1}+\mu_r\right) & \mu_{r} &\cdots  \\
\cdots & \cdots & \cdots & \cdots & \cdots  & \cdots &\cdots \\
\end{array}
\right).} \label{class1-inf}
\end{equation}
\noindent In the case of finite state space (i.e. $S<\infty$) the matrix $B^*(t)=$
\begin{equation}
{\scriptsize \left(
\begin{array}{ccccc}
-\left(\lambda_0+\mu_1\right)  & \mu_1
 & 0 & \cdots & 0 \\
\lambda_1  & -\left(\lambda_1+\mu_2\right) & \mu_2 & \cdots & 0 \\
\ddots & \ddots & \ddots & \ddots & \ddots  \\
0 & \cdots & \cdots & \lambda_{S-1} &
-\left(\lambda_{S-1}+\mu_S\right)
\end{array}
\right).} \label{class1-fin}
\end{equation}
\noindent Note that the matrix $B^*(t)$ is essentially
non-negative for any $t \ge 0$ i.e.
all its off-diagonal elements are non-negative for any $t$.

\subsection{$B^*(t)$ for $X(t)$ belonging to class~$\rm (ii)$}

Consider a process $X(t)$ with $a_{ij}(t)=0$ for $i<j-1$,
$a_{i+k,i}(t)=a_k(t)$ for $k\ge 1$ and $a_{i,i+1}(t)=\mu_{i+1}(t)$.
Such $X(t)$ describes the evolution of the total number of customers
in a queue with batch arrivals and single services ($a_k(t)$ are the
(state-independent) intensities of group arrivals and $\mu_{i+1}(t)$
are the (state-dependent) service intensities). Such processes in
the simplest forms were firstly considered in \cite{Nelson1988} and,
under the assumption of decreasing of $a_k(t)$, have been studied in
\cite{Zeifman2018c}. In the case of countable state space (i.e.
$S=\infty$) the matrix $B^*(t)$ has the form:

\begin{equation}
{\scriptsize B^*(t) =  \left(
\begin{array}{cccccc}
a_{11}  & \mu_1
 & 0 & \cdots & 0 \\
a_1  & a_{22} & \mu_2 & \cdots & 0 \\
a_2  & a_1  &a_{33} & \mu_3 & \cdots &  \\
\ddots & \ddots & \ddots & \ddots & \ddots  \\
\ddots & \ddots & \ddots & \ddots & \ddots  \\
\end{array}
\right)}. \label{class2-inf}
\end{equation}
\noindent
 In the case of finite state space (i.e. $S<\infty$) the matrix $B^*(t)=$
\begin{equation}
{\scriptsize \left(
\begin{array}{ccccc}
a_{11}-a_S  & \mu_1
 & 0 & \cdots & 0 \\
a_1-a_S  & a_{22}-a_{S-1} & \mu_2 & \cdots & 0 \\
\ddots & \ddots & \ddots & \ddots & \ddots  \\
a_{S-1}-a_S & \cdots & \cdots & a_1-a_2 &
a_{SS}-a_1
\end{array}
\right)}. \label{class2-fin}
\end{equation}

\noindent Note that the matrix $B^*(t)$ is essentially
non-negative for any $t \ge 0$ if the arrival intensities $a_k(t)$ are
decreasing in $k$.

\subsection{$B^*(t)$ for $X(t)$ belonging to class~$\rm (iii)$}

Consider a process $X(t)$ with
$a_{ij}(t)=0$ for $i>j+1$, $a_{i,i+k}(t)=b_k(t)$, $k\ge
1$ and $a_{i+1,i}(t)=\la_i(t)$.
Such $X(t)$ describes the evolution
of the total number of customers
in a queue with batch services and single arrivals
($\la_i(t)$ are the (state-independent) arrival intensities
and $b_k(t)$ are the (state-independent) intensities
of service of a group of $k$ customers).
Such processes were considered to some extent in
\cite{Nelson1988,Li2017}.
In the case of countable state space (i.e. $S=\infty$)
the matrix $B^*(t)=$
\begin{equation}
{\tiny \left(
\begin{array}{cccccc}
-\left(\lambda_0 +b_1\right)  & b_1 - b_2
 & b_2 - b_3 & \cdots & \cdots \\
\lambda_1 & -\big(\lambda_1+\sum\limits_{i\le 2}b_i\big) & b_1 - b_3 & \cdots & \cdots\\
\ddots & \ddots & \ddots & \ddots & \ddots  \\
0 & \cdots & \cdots & \lambda_{r-1} &
-\big(\lambda_{r-1}+\sum\limits_{i\le r}b_i\big) \cdots \\
\ddots & \ddots & \ddots & \ddots & \ddots  \\
\end{array}
\right)}.\label{class3-inf}
\end{equation}
 \noindent  In the case of finite state space (i.e. $S<\infty$)
the matrix $B^*(t)$ is
$$B^*(t) =$$
\begin{equation}
{\tiny  \left(
\begin{array}{ccccc}
-\left(\lambda_0 +b_1\right)  & b_1 - b_2
 & b_2 - b_3 & \cdots & b_{S-1} - b_S \\
\lambda_1  & -\big(\lambda_1+\sum\limits_{i\le 2}b_i\big) & b_1 - b_3 & \cdots & b_{S-2} - b_S \\
\ddots & \ddots & \ddots & \ddots & \ddots  \\
0 & \cdots & \cdots & \lambda_{S-1} &
-\big(\lambda_{S-1}+\sum\limits_{i\le S}b_i\big)
\end{array}
\right)}.\label{class3-fin}
\end{equation}

\noindent Note that the matrix $B^*(t)$ is essentially
non-negative for any $t \ge 0$ if the service intensities $b_k(t)$ are
decreasing in $k$.

\subsection{$B^*(t)$ for $X(t)$ belonging to class~$\rm (iv)$}

Consider a process $X(t)$ with $a_{i+k,i}(t)=a_k(t)$ and
$a_{i,i+k}(t)=b_k(t)$ for ${k\ge 1}$. Such $X(t)$ describes the
evolution of the total number of customers in an inhomogeneous queue
with (state-independent) batch arrivals and group services ($a_k(t)$
are the (state-independent) intensities of group arrivals and
$b_k(t)$ are the (state-independent) intensities of group services).
Such process under the assumption of decreasing in $k$ intensities
$a_k(t)$ and $b_k(t)$ have been studied in \cite{Zeifman2014b}. In
the case of countable state space (i.e. $S=\infty$) the matrix
$B^*(t)$ has the form:

\begin{equation}
B^* = {\small \left(
\begin{array}{cccccc}
a_{11}  &  b_1 - b_2
 & b_2 - b_3 & \cdots & \cdots \\
a_1  & a_{22} & b_1 - b_3 & \cdots & \cdots \\ \\
\ddots & \ddots & \ddots & \ddots & \ddots  \\
a_{r-1} & \cdots & \cdots & a_1 & a_{rr} & \cdots \\
 \cdots & \cdots & \cdots & \cdots & \cdots & \cdots \\
\end{array}
\right)}.\label{class4-inf}
\end{equation}
 \noindent In the case of finite state space (i.e. $S<\infty$)
the matrix $B^*(t)=$

\begin{equation}
{\small \left(
\begin{array}{ccccc}
a_{11}-a_S  &  b_1 - b_2
 & b_2 - b_3 & \cdots & b_{S-1} - b_S\\
a_1-a_S  & a_{22}-a_{S-1} & b_1 - b_3 & \cdots &
b_{S-2} - b_S \\
\ddots & \ddots & \ddots & \ddots & \ddots  \\
a_{S-1}-a_S & \cdots & \cdots & a_1-a_2 &
a_{SS}-a_1
\end{array}
\right)}.\label{class4-fin}
\end{equation}
\noindent Note that the matrix $B^*(t)$ is essentially
non-negative for any $t \ge 0$ if the intensities $a_k(t)$ and $b_k(t)$ are
decreasing in $k$.

\section{Upper bounds using logarithmic norm}

Throughout this section by $\|\cdot\|$ we denote  the $l_1$-norm,
i.e. $\|{\vp(t)}\|=\sum_{i\in \mathcal{X}} |p_i(t)|$ and $\|A(t)\| =
\sup_{j \in \mathcal{X}} \sum_{i\in \mathcal{X}} |a_{ij}(t)|$. Let
$\Omega$ be a set of all stochastic vectors, i.e. $l_1$ vectors with
non-negative coordinates and unit norm. Recall that a Markov chain
$X(t)$ is called {\it weakly ergodic}, if
${\|\vp^{*}(t)-\vp^{**}(t)\| \to 0}$ as $t \to \infty$ for any
initial conditions $\vp^{*}(0)$ and $\vp^{**}(0)$, where
$\vp^{*}(t)$ and $\vp^{**}(t)$ are the corresponding solutions of
(\ref{ur01}).

Recall that the logarithmic norm\footnote{A number of queueing
applications of this approach have been studied in
\cite{Granovsky2004,Zeifman2006,Zeifman2018c}.} of the operator
function ${B}(t)$ is defined as
$$
\gamma({B}(t)) = \lim_{h \to
+0}h^{-1}\left(\|I+hB(t)\|-1\right).
$$

\noindent Denote by $V(t, s)= V(t)V^{-1}(s)$ the Cauchy operator of the equation (\ref{hom1}).
Then $\|V(t, s)\| \le e^{\, \int_s^{t} \gamma(B(u))\, du}$. For an operator function from $l_1$
to itself we have the formula
\begin{equation}
\gamma({B}(t)) = \sup_{j\in \mathcal{X}} \left(b_{jj}(t)+\sum_{i \in \mathcal{X}, i \neq j} |b_{ij}(t)|\right).
\label{lognorm1}
\end{equation}

\noindent {Note that if the matrix $B(t)$ is essentially
non-negative then $\gamma({B}(t)) = \sup_{j \in \mathcal{X}} \left
(\sum_{i \in \mathcal{X} } b_{ij}(t) \right )$.}

%\noindent Moreover, if the matrix $B(t)$ is essentially non-negative then
%\begin{equation}
%\gamma({B}(t)) = \sup_j \left(b_{jj}(t)+\sum_{i \neq j}
%b_{ij}(t)\right). \label{lognorm2}
%\end{equation}

Assume that the state space $\mathcal{X}$ is countable i.e.
$S=\infty$. Let ${\{d_i, \ i \ge 1\}}$ be a sequence of positive
numbers and let $D=diag(d_1,d_2,\dots)$ be the diagonal matrix, with
the off-diagonal elements equal to zero. By putting ${\bf w}(t)=D
{\bf u}(t)$ in (\ref{hom11}), we obtain the following equation
\begin{equation}
\frac{d }{dt} {\bf w}(t)= B^{**}(t){\bf w}(t), \label{hom111}
\end{equation}
\noindent where $B^{**}(t)=D
B(t)^{*}D^{-1}$.
%=\left(b^{**}_{ij}(t)\right)_{i,j=1}^S$.
Put\footnote{It is possible to obtain explicit expressions for
$\alpha_i(t)$ for all of the considered classes (i)--(iv) (see the
details in \cite{Zeifman2018c}).}
\begin{equation}
\alpha_i\left(t\right)= -\sum_{j=1}^\infty b_{ji}^{**}(t),
 \ i \ge 1, \label{posit02}
\end{equation}
\noindent and let $\alpha(t)$ and $\beta(t)$ denote the
least lower and the least upper bound of the
sequence of functions $\{ \alpha_i(t), \ i \ge 1\}$ i.e.
\begin{equation}
\alpha\left(t\right)=\inf_{i \ge 1}\alpha_i\left(t\right), \ \beta
\left(t\right)=\sup_{i \ge 1}\alpha_i\left(t\right). \label{posit03}
\end{equation}

\noindent The next theorem and corollary have been proved in \cite[Theorem~1]{Zeifman2018c}
and are stated here for the sake of completeness.

\begin{theorem}{}
Let there exist a sequence  $\{d_i, \ i \ge 1\}$ of positive numbers
such that $d_1=1$, $d=\inf_{i \ge 1} d_i > 0$ and $B^*(t)$ is essentially
non-negative. Let $\alpha(t)$, defined by (\ref{posit03}), satisfy
\begin{equation}
\int_0^{\infty} \alpha(t)\, dt = + \infty.
\end{equation}

\noindent Then the Markov chain $X(t)$ is weakly ergodic and
for any initial conditions $s \ge 0$, ${\bf w}(s)$ and any  $t \ge s$
the following upper bound holds:
\begin{eqnarray}
\|{\bf w}\left(t\right)\| \le e^{-\int_s^t
{\alpha\left(u\right)du}}\|{\bf w}(s)\|. \label{t001}
\end{eqnarray}
\noindent
If in addition all components of the vector ${\bf w}\left(s\right)$ are non-negative, then
for any $0 \le s \le t$ the following lower bound holds:
\begin{equation}
\|{\bf w}\left(t\right)\| \ge
e^{-\int_s^t {\beta\left(u\right)du}}\|{\bf w}(s)\|. \label{t004}
\end{equation}
\end{theorem}

\medskip

\begin{corollary}{}%\hspace{-0.2cm}{\bf }
Let under the assumptions of Theorem 1 the sequence ${\{d_i, \ i \ge
1\}}$ be such that all $\alpha_i\left(t\right)$ do not depend on
$i$ i.e. are the same for any $i$. Then
$\alpha\left(t\right)= \beta\left(t\right)$ and the upper bound
(\ref{t001}) is sharp. If in addition all components of the vector
${\bf w}(s)$ are non-negative, for any  $0 \le s\le t$ it holds that
\begin{eqnarray}
\|{\bf w}\left(t\right)\| = e^{-\int_s^t
{\alpha\left(u\right)du}}\|{\bf w}\left(s\right)\|. \label{t001eq}
\end{eqnarray}

\end{corollary}

\medskip

If the Markov chain $X(t)$
is homogeneous, then the expressions in (\ref{posit02})
and (\ref{posit03}) do not depend on $t$.
In such case the upper  and lower bounds \eqref{t001}, \eqref{t004} can be improved.
The following result is due to \cite[Theorem~2]{Zeifman2018c}.

\begin{theorem}{}%\hspace{-0.2cm}{\bf }
Let there exist a sequence $\{d_i, \ i \ge 1\}$ of positive numbers
such that $d_1=1$, $d=\inf_{i\ge 1} d_i > 0$ and $B^*(t)$ is essentially
non-negative. Let $\alpha$, defined by \eqref{posit03}, be positive.
Then $X(t)$ is ergodic and for any initial condition ${\bf w}(0)$ and any $t \ge 0$
 the following upper bound holds:
\begin{eqnarray}
\|{\bf w}\left(t\right)\| \le e^{- \alpha t}\|{\bf
w}\left(0\right)\|. \label{t001hom}
\end{eqnarray}
\noindent If in addition all components of the vector ${\bf w}\left(0\right)$ are non-negative, then
for any $ t \ge 0$ the following lower bound holds:
\begin{equation}
\|{\bf w}\left(t\right)\| \ge e^{- \beta t}\|{\bf
w}\left(0\right)\|. \label{t004hom}
\end{equation}
\noindent If $\alpha=\beta$ then the bound (\ref{t001hom}) is sharp.
\end{theorem}

Assume now that the state space is finite i.e. $S <\infty$. Then
$d_i$ can be arbitrary positive numbers and we can find such
constants, say $C_1$ and $C_2$, that

$$
\|{\bf w}(t)\| = \|DT{\bf y}(t)\| \le C_1\|{\bf y}(t)\|,$$ $$
 \ \|{\bf y}(t)\| =
\|T^{-1}D^{-1}{\bf w}(t)\| \le C_2\|{\bf w}(t)\|. %\label{norms1}
$$
\noindent
%where $\|\cdot\|$ denotes $l_1$ norm which is closely
%related with total variance distance.
Hence {\it Theorem}~1 and {\it Theorem}~2 provide bounds on the rate
of convergence in the $l_1$-norm. The explicit expressions for the
constants can be found in \cite{Granovsky2004,Zeifman2006}. If the
Markov chain $X(t)$ is homogeneous and $\alpha^*$ is the decay
parameter, defined as
$$
\lim_{t \rightarrow \infty} (p_{ij}(t)-\pi_j)= O(e^{-\alpha^* t}),
$$
\noindent
where $\{ \pi_j, \ j \ge 0\}$ are the stationary probabilities of
the chain, then ${\alpha \le \alpha^* \le \beta}$.

Notice that some additional results for finite homogeneous Markov
chains $X(t)$ belonging to class (i) we can find in \cite{VanDoorn2010}. In
particular, in \cite{VanDoorn2010} it was proved that the exact estimate of the
rate of convergence can be obtained. In the next theorem we provide an alternative
proof of this fact.

\begin{theorem}{}%\hspace{-0.2cm}{\bf }
Let $X(t)$ be a homogeneous birth-death processes
with a finite state space of size $S$ and
let all birth and death intensities be positive.
Then there exists a set $\{d_i, \ 1\le  i \le S\}$
of positive numbers such that  $\alpha = \alpha^* = \beta$,
where $\alpha^*$ is the decay parameter of $X(t)$,
and $\alpha$ and $\beta$ are defined by \eqref{posit03}.
\end{theorem}

{\bf Proof.}
% This follows from the next statement.
%\medskip
%\begin{lemma}\hspace{-0.2cm}{\bf .}
Let $C$ be an essentially non-negative irreducible matrix
such that there exists $n_0>0$ such that $C^{n_0} >0$.
%Let some degree of $C$ be positive matrix.
Denote by $\lambda_0$ its maximal eigenvalue. It is simple and positive.
Then there exists a diagonal matrix with positive entries $D= diag (d_1,
\dots, d_S)$ such that all column sums for matrix $C_{D}=DCD^{-1}$
are equal to $\lambda_0$.
%\end{lemma}
%\medskip
%{\bf Proof.}
Indeed, let $m= \max_{1 \le j \le S}{|c_{jj}|}$. Consider the irreducible matrix
$C'=C^{T}+mI$. It has a simple eigenvalue  $\lambda^*=\lambda_0+m$
and the corresponding eigenvector ${\bf x} = \left(x_{1}, \dots,
x_{S}\right)^T$ has strictly  positive coordinates.
Put $d_{i}=x_{i}^{-1}$, $1\le  i \le S$.
Then  ${\bf e} = \left(1, \dots, 1\right)^T$ is
the eigenvector of the matrix  $C'_{D}=DC'D^{-1}$.
Therefore all row sums in the matrix $C'_{D}$ are equal to $\lambda^*$.
Thus all row sums in the matrix  $C^{T}_{D}= C'_{D} -mI$ are equal to $\lambda^* -m=\lambda_0$,
and all column sums of the matrix $C_{D}$ are equal to $\lambda_0$.

%\begin{flushright}
%\qedsymbol
%\end{flushright}

\section{Upper bounds using Lyapunov functions}

As is was mentioned in the introduction the method based on Lyapunov functions
no longer applies in the case of countable state space $\mathcal{X}$.
In this section, under the assumption that $\mathcal{X}$ is finite
i.e. $S< \infty$, it is shown how (quadratic) Lyapunov functions
can be applied to obtain the explicit upper bounds on the rate of convergence
of some $X(t)$ belonging to classes (i)--(iii).
Unlike the bounds provided by the method based on the logarithmic norm,
Lyapunov functions yield bounds in $l_2$-norm (Euclidean norm) and
thus are somewhat weaker.

Denote throughout this section by $\|\cdot\|$
the $l_2$-norm, i. e. $\|{\vp}(t)\|=\sqrt{\sum_{i \in \mathcal{X}}p_i(t)^2}$.
Consider the system (\ref{hom111}). Let $V(t)=\sum_{k=1}^S w^2_k (t)$, where
${\bf w}(t) = \left( {w}_1(t), w_2(t), \dots, w_S(t)\right)^T$
is the solution of (\ref{hom111}). By differentiating $V(t)$
we obtain
\begin{eqnarray}
\frac{dV(t)}{dt}= \sum_{k=1}^S 2w_k(t)\frac{dw_k(t)}{dt}= \\
=-2\sum_{i=1}^S \sum_{j=1}^S \left(-b_{ij}^{**}(t)\right)w_i (t) w_j(t). \nonumber
 \label{class1**1}
\end{eqnarray}

\noindent If we find a set of positive numbers $\{d_i, 1 \le i \le S\}$ and
a function $\beta^*(t)$ satisfying
\begin{eqnarray}
\frac{dV(t)}{dt} \le -2\beta^*(t)V(t),
 \label{lyap02}
\end{eqnarray}
\noindent for any ${\bf w}(t)$, being the solution of (\ref{hom111}), then
for a $X(t)$ belonging to classes (i)--(iv)
and for any initial condition ${\bf w}(0)$ it holds that
\begin{equation}
\|{\bf w}(t)\| \le e^{- \int_s^t\beta^*(\tau)\,d{\tau}}\|{\bf
w}(0)\|.
 \label{lyap03}
\end{equation}

\noindent
For a finite homogeneous Markov chain $X(t)$ belonging to class (i)
such set $\{d_i, 1 \le i \le S\}$ is given in the next theorem.

\begin{theorem}{}%\hspace{-0.2cm}{\bf}
Let $X(t)$ be a homogeneous birth-death process
defined on a finite state space $\mathcal{X}$
with the possibly state-dependent birth intensities $\lambda_{k}$
and possibly state-dependent death intensities $\mu_{k}$.
Assume that $\lambda_{k}>0$ and $\mu_{k}>0$ for each $k \in \mathcal{X}$.
Then there exist a set of positive numbers
$\{d_i, 1 \le i \le S\}$, a positive number $\beta^{*}$
and a set of numbers $\{ \alpha_i, 1 \le i \le S\}$ such that
\begin{equation}
\frac{dV(t)}{dt} = -2 \beta^{*} \sum^{S}_{k=1} w^{2}_k -
2 \sum^{S-1}_{k=1} (\alpha_k  w_k - \alpha_{k+1}
w_{k+1})^2.
 \label{lyap11}
 \end{equation}
\end{theorem}

\noindent
{\bf Proof.} If $X(t)$ is a homogeneous birth-death process, then $B^*(t)$
does not depend on $t$ and thus it is constant tridiagonal matrix.
Let $d_1=1$, ${d_{k+1}=d_k\sqrt{\mu_k /\la_k}}$, $k \ge 1$.
Remembering that $D = diag(d_1,\dots,d_S)$
and $B^{**}(t)=D B(t)^{*}D^{-1}$, we immediately obtain $B^{**}=$
\begin{equation}
{\tiny\left(
\begin{array}{ccccc}
-\left(\lambda_0+\mu_1\right)  & \sqrt{\la_1\mu_1}
 & 0 & \cdots & 0 \\
\sqrt{\la_1\mu_1}  & -\left(\lambda_1+\mu_2\right) & \sqrt{\la_2\mu_2} & \cdots & 0 \\
\ddots & \ddots & \ddots & \ddots & \ddots  \\
0 & \cdots & \cdots & \sqrt{\la_{S-1}\mu_{S-1}} &
-\left(\lambda_{S-1}+\mu_S\right)
\end{array}\right)}.
\end{equation}
\noindent Note that $B^{**}$ is the symmetric matrix.
By putting $\Phi(t)= -0.5 \frac{dV(t)}{dt}$ in \eqref{class1**1}
we obtain
$$
\Phi(t)=\lambda_0
w_1^2+\mu_S  w_S^2 + \sum^{S-1}_{k=1} (\sqrt{\mu_k}  w_k -
\sqrt{\lambda_k}  w_{k+1})^2.
$$

\noindent Choose a positive number $\beta$
such that $\beta<\min(\lambda_0,\lambda_1,\dots,\lambda_S)$
and put ${\phi_0=\lambda_0-\beta}$. Then the terms in the right-hand side
of the previous relation can be rearranged to give
$$
\Phi(t)=\beta
w_1^2 + \left(\sqrt{\mu_1+\phi_0}
 w_1 -
\frac{\sqrt{\lambda_1 \mu_1}}{\sqrt{\mu_1+\phi_0}} w_2 \right)^2+$$
$$+\lambda_1 \left(\frac{\phi_0}{\mu_1+\phi_0}\right) w_2^2+$$
$$+\sum^{S-1}_{k=2} (\sqrt{\mu_k}  w_k - \sqrt{\lambda_k}  w_{k+1})^2 +\mu_S  w_S^2.
$$

\noindent Consider the coefficient of $w_2^2$.
Note that it can always\footnote{Indeed, if $\beta$ is larger than
the coefficient of $w_2^2$, it suffices to make one step back
and choose a new value of $\beta$
(satisfying $\beta<\min(\lambda_0,\lambda_1,\dots,\lambda_S)$)
smaller than the current one.} be represented as $\beta+\phi_1$
with $\phi_1>0$. Thus we can rearrange the terms
in the previous relation and obtain
$$
\Phi(t) =\beta (w_1^2+w_2^2)+ \left(\sqrt{\mu_1+\phi_0} w_1 - \frac{\sqrt{\lambda_1 \mu_1}}{\sqrt{\mu_1+\phi_0}}  w_2 \right)^2+
$$
$$+\left(\sqrt{\mu_2+\phi_1}  w_2 - \frac{\sqrt{\lambda_2 \mu_2}}{\sqrt{\mu_2+\phi_1}} w_3 \right)^2
+\lambda_2 \left(\frac{\phi_1}{\mu_2+\phi_1}\right)
w_3^2+$$
$$+\sum^{S-1}_{k=3} (\sqrt{\mu_k}  w_k - \sqrt{\lambda_k}  w_{k+1})^2 +\mu_S  w_S^2.$$

\noindent By proceeding in the similar manner (i.e. choosing suitable value of $\beta$,
representing each coefficient of $w_k$ as $\beta+\phi_{k-1}$, $\phi_{k-1}>0$, and rearranging the terms)
we arrive at the following representation of $\Phi(t)$:
$$
\Phi(t)=\beta  \sum^{S-1}_{k=1} w_k^2+$$ $$\sum^{S-1}_{k=1}\left(\sqrt{\mu_k+\phi_{k-1}}
 w_k -
\frac{\sqrt{\lambda_k \mu_k}}{\sqrt{\mu_k+\phi_{k-1}}} w_{k+1}
\right)^2+$$
$$+\left(\mu_S +\lambda_{S-1} \frac{\phi_{S-2}}{\mu_{S-1}+\phi_{S-2}}\right)  w_S^2.
$$

\noindent
If the coefficient of $w_S^2$ is larger than $\beta$, then
we can choose any $\beta^\ast$ such that
$\beta^\ast \in \left(\beta;\mu_S +\lambda_{S-1}\phi_{S-2}/(\mu_{S-1}+\phi_{S-2})\right)$.
Therefore $\beta<\beta^\ast<\beta+\varepsilon$, where
$\beta+\varepsilon =\mu_S +\lambda_{S-1}\phi_{S-2}/(\mu_{S-1}+\phi_{S-2})$. Put
$\beta_{1}= \beta+ 0.5\varepsilon$ and continue the process of
selecting squares in the opposite direction (starting from $w_S^2$).
If the coefficient of $w_S^2$ is less than $\beta$, then we can choose
$\beta^\ast \in \left(\beta-\varepsilon; \beta \right)$. In this
case we put $\beta_{1}= \beta-0.5\varepsilon$ and continue the
process of selecting squares, starting from $w_1^2$.
In such a way we get a sequence of nested segments converging to $\beta^\ast$.

%\begin{flushright}
%\qedsymbol
%\end{flushright}

Note that the existence of the upper bound
$\|{\bf w}(t)\| \le e^{- \beta^* t}\|{\bf w}(0)\|$
also follows from \eqref{lyap02}, \eqref{lyap03} and
the Theorem~4. The inequality turns into the equality
once the set of numbers $\{ \alpha_i, 1 \le i \le S\}$
is chosen in such a way that the second sum in \eqref{lyap11} is equal to zero.

Let us specify the upper bound \eqref{lyap03}
for some finite homogeneous Markov chains $X(t)$ belonging to class (ii).
Specifically let in a process $X(t)$ belonging to (ii)
the arrival intensities be such that
$\la_1=0$ and $\la_k=\la$ for $2 \le k \le S$.
From the queueing perspective this means that only
arrivals in batches are possible.
Then the matrix $B^*(t)$ given by (\ref{class2-fin})
does not depend on $t$ and takes the following form:
\begin{equation} { B^* =  \left(
\begin{array}{ccccccc}
a_{11}-\la  &  \mu_1   & 0 & \cdots &  0 \\
-\la  &  a_{22}-\la  &\mu_2  & \cdots &  0 \\
0  &  -\la  &a_{33}-\la   & \cdots &  0  \\
\cdots \\
0  & 0  & 0  & \cdots & a_{SS} \\
\end{array}
\right)}. \label{class2*}
\end{equation}

\noindent Let $d_1=1$, $d_{k+1}= d_k \sqrt{\mu_k/\la}$, $k \ge 1$.
Remembering that $D = diag(d_1,\dots,d_S)$
and $B^{**}(t)=D B(t)^{*}D^{-1}$, we immediately obtain $B^{**} =$
{\scriptsize
\begin{eqnarray}
 { \left(
\begin{array}{ccccc}
-\left(S\la+\mu_1\right)  & \sqrt{\la\mu_1}
 & 0 & \cdots & 0 \\
-\sqrt{\la\mu_1}  & -\left((S-1)\la+\mu_2\right) & \sqrt{\la\mu_2} & \cdots & 0 \\
\ddots & \ddots & \ddots & \ddots & \ddots  \\
0 & \cdots & \cdots & -\sqrt{\la\mu_{S-1}} & -\mu_S
\end{array}
\right)}.\nonumber \label{class2**}
\end{eqnarray}}

\noindent For such matrix $B^{**}$ the equation \eqref{class1**1}
can be rewritten as $\frac{dV(t)}{dt}=2 \sum^{S}_{k=1} (a_{kk}-\la)w_k^2(t)$,
%\begin{eqnarray}
%\frac{dV(t)}{dt}=
%\rred{2 \sum^{S}_{k=1} (a_{kk}-\la)w_k^2(t),}
%2\left((a_{11}-\la)w_1^2(t) +(a_{22}-\la)w_2^2(t) + \dots +a_{SS}w_S^2(t)\right),
% \label{class2**1}
%\end{eqnarray}
%\noindent
wherefrom the next theorem follows\footnote{Note that in the
considered case we can also obtain the lower bound on the rate of
convergence using the approach in \cite{Zeifman2018spl}}.

\begin{theorem}{}%\hspace{-0.2cm}{\bf }
Let $X(t)$ be a homogeneous Markov chain
defined on a finite state space $\mathcal{X}$
with the state-independent group arrival intensities $q_{k,k+i}=\lambda$, $i\ge 2$, $q_{k,k+1}=0$,
and possibly state-dependent service intensities $q_{k,k-1}=\mu_{k}$, $1 \le k \le S$.
Then the following bound on the rate of convergence holds:
\begin{equation}
 \|{\bf w}(t)\| \le e^{-
\beta^*t}\|{\bf w}(0)\|,
 \label{class2**3}
\end{equation}
\noindent where $\beta^*=\min\left(S\la+\mu_1,\dots, 2\la+\mu_{S-1},\mu_{S}\right)$
%\begin{eqnarray}
%,
% \label{class2**4}
%\end{eqnarray}
i.e. $\beta^*$ is the decay parameter (spectral gap) of the Markov chain.
\end{theorem}

% далее не вносил изменений

\noindent
Note that a similar result can be obtained for
the homogeneous Markov chains $X(t)$ belonging to class (iii).
The following example shows that
Lyapunov functions lead to explicit uppers bounds for the rate of convergence
also for finite inhomogeneous Markov chains.

\bigskip

\noindent
{\bf Example 1.} Consider the Markov process $X(t)$ that describes
the evolution of the total number of customers in the $M(t)/M(t)/1/S$ queue with bulk arrivals,
when all transition intensities are periodic functions of time.
Let the arrival intensities
be $a_1(t)=1+\sin 2\pi t$, $a_k(t)=2+\sin 2\pi t+\cos 2\pi t$ for $2 \le k \le S$ and all
the service intensities be $\mu_k(t)=m^2\left(1+\cos 2\pi t\right)$ for $1 \le k \le S$ and some $m \ge 1$.
Such $X(t)$ belongs to class (ii).
By setting $d_1=1$, $d_{k+1}= m d_k $, $k \ge 1$, we obtain the matrix  $B^{**}(t)$
in the following form $B^{**}(t)=(b_{ij}^{**}(t))$, where
$$b_{i,i+1}^{**}(t)=m\left(1+\cos 2\pi t\right),$$
$$b_{i,i}^{**}(t)=a_{ii}(t)-a_{S-i+1},$$
$$b_{i+1,i}^{**}(t)=-m\left(1+\cos 2\pi t\right).$$
%$$b_{i,S}^{**}(t)=-m^{-2}\left(2+\cos(2\pi t)\right) \frac{d_{i}}{d_{S}}, ~~  i< S,$$
%$$b_{SS}^{**}(t)=-10(2+\sin(2\pi t)) - m^{-2}\left(2+\cos(2\pi t)\right), ~~ i=S.$$

\noindent Then $\frac{dV(t)}{dt}=2 \sum^{S}_{k=1} (a_{kk}(t)-a_{S+1-k}(t))w_k^2(t)$
and from (\ref{class2**}) it follows that
%\begin{equation*}
%\frac{dV(t)}{dt}= 2\left((a_{11}(t)-a_S(t))w_1^2(t)
%+(a_{22}(t)-a_{S-1}(t))w_2^2 (t) + \dots
%+(a_{SS}(t)-a_1(t))w_S^2(t)\right),
% \label{lyap22}
%\end{equation*}
%\noindent
for any initial condition ${\bf w}(0)$ the sharp upper bound on the rate of convergence is
\begin{equation*}
\|{\bf w}(t)\| \le e^{- \int_0^t\beta^*(\tau)\,d{\tau}}\|{\bf
w}(0)\|,
 \label{lyap03new}
\end{equation*}
\noindent where $\beta^*(t)=2+\sin 2\pi t+\cos 2\pi t$. Note that,
for the considered case, the method based on Lyapunov functions
yields the best (among the three methods considered in this paper)
possible upper bound. It is also worth noticing that we can apply
the obtained upper bound for the computation of the limiting
distribution of $X(t)$. For example, let $S=199$ and $m=90$. Then,
using truncation techniques, which were developed in
\cite{Zeifman2006,Zeifman2014c}, any limiting probability
characteristic of $X(t)$ can be computed  with the given
approximation error. In Figs.~1--8 we can see the behaviour of the
conditional expected number $E(X(t)|X(0))$ of customers in the queue
at instant $t$ and the state probabilities $p_0(t)$, $p_{99}(t)$ and
$p_{199}(t)$ as functions of time $t$ under different initial
conditions $X(0)$. The approximation error is~$10^{-3}$.

Note that one general framework for the computation of the
limiting characteristics  of time-dependent queueing systems is
described in detail in the recent paper \cite{Satin2019}.
Particularly, having the bounds on the rate of convergence we can
 compute the time instant, say $t^*$, starting from which
probabilistic properties of $X(t)$ do not depend on the value of
$X(0)$ (assuming that  the process starts at time $t=0$). Thus, for
example, if the transition intensities are periodic (say,
1-time-periodic), we can truncate the process on the interval
$[t^*, t^*+1]$ and solve the forward Kolmogorov system of
differential equations on this interval with $X(0)=0$. In such a
way, we can build approximations for any limiting probability
characteristics of $X(t)$ and estimate stability (perturbation)
bounds.

%an approximation of
%important characteristics of the queue-length process, such as the
%mean (the mathematical expectation $E(t,k)$)  under different
%initial conditions and some probabilities of different fixed queue
%lengths. Firstly, on the final interval $[t^*, t^*+1]$ and secondly,
%on interval $[0, t^*]$. It is interesting to note the strange
%behavior of state probabilities.

\section{Upper bounds using differential inequalities}

As it was firstly shown in \cite{Zeifman2019rodos},
there are situations when the previous two methods for bounding
the rate of convergence do not work well
(either lead to poor upper bounds or do not yield
upper bounds at all).
Here we present probably the most general method,
which is based on differential inequalities
and which can be applied to a $X(t)$ belonging to classes (i)--(iv)
with finite state space (i.e. $S< \infty$)
and all transition intensity functions being analytic  functions of time $t$.

Throughout this section by $\|\cdot\|$ we denote the $l_1$-norm.
Consider a finite system of linear differential equations
\begin{equation}
 \frac{d}{dt}\vx(t)=A(t)\vx(t), \ t \ge 0,
\label{01}
\end{equation}
\noindent where $A(t)$ is some matrix\footnote{This matrix $A(t)$
must not be confused with the matrix in \eqref{ur01}.} with all entries $a_{ij}(t)$ being analytic functions of $t$ and ${ \vx}(t) = \left({x}_1(t), \dots, {x}_S(t)\right)^T$.
Let $\vx(t)$ be an arbitrary solution of (\ref{01}). Consider an
interval $[t_1,t_2]$ with fixed signs of coordinates of $\vx(t)$
(i.e. $x_i(t)\neq 0$ for all $1 \le i \le S$ and for all $t\in [t_1,t_2]$).
Choose the set of numbers $\{d_i, 1 \le i \le S\}$
such that the sign of each $d_i$ coincides with the sign of $x_i(t)$.
Then $d_ix_i(t) \ge 0$ for all $t \in [t_1,t_2]$ and
hence $\sum_{k=1}^{S} d_k x_k (t) = \|\vx(t)\|$ can be considered as the
$l_1$-norm.

Put $\vz(t)=D\vx(t)$ and $\tilde{A}(t)=DA(t)D^{-1}$,
%=\left(\tilde{a}_{ij}(t)\right)_{i,j=1}^S$,
where $D = diag(d_1,\dots,d_S)$,
and consider the following system of differential equations
\begin{equation}
 \frac{d}{dt}\vz(t)=\tilde{A}(t)\vz(t),
\label{02}
\end{equation}
\noindent for  $t \in [t_1,t_2]$.
%\begin{equation}
%\sum_i\tilde{a}_{ij}(t) \le  -\alpha_D(t), \quad j=1,\dots, S.
%\label{03}
%\end{equation}
If for the chosen matrix $D$ there exists
a function\footnote{The lower index in $\alpha_D(t)$ is used to explicitly indicate
that this function depends on choice of the matrix $D$.} $\alpha_D(t)$ such that
$\sum_{i=1}^S\tilde{a}_{ij}(t) \le -\alpha_D(t)$ for each $1 \le j \le S$,
then the following bound holds
\begin{equation}
\frac{d}{dt}\|\vz(t)\|=
%\frac{d\left(\sum_kz_k\right)}{dt}=
\sum_{j=1}^S\sum_{i=1}^S \tilde{a}_{ij}(t)z_j(t)
\le -\alpha_D(t)\|\vz(t)\|. \label{04}
\end{equation}

\noindent Choose $\alpha^*(t)$ such that $\alpha^*(t)= \min\alpha_D(t)$,
where the minimum is taken
over all time intervals $[t_1,t_2]$, $0<t_1<t_2$, with different combinations of coordinate signs
of the solution ${ \vx}(t)$. For any such combination there exists
a particular inequality $\|\vz(t)\| \le e^{-\int_{t_1}^{t_2}\alpha^*(\tau)\,
d\tau}\|\vz(t_1)\|$.

%Let us now compare all the norms.
From the fact that there exist constants, say $C_1$ and $C_2$,
such that  $\|\vx(t)\| \le C_1 \|\vz(t)\|$
and $\|\vz(t)\| \le C_2\|\vx(t)\|$ for
any interval $[t_1,t_2]$, $0<t_1<t_2$, and any corresponding diagonal matrix
$D$, the following theorem follows.

\begin{theorem}{}%\hspace{-0.2cm}{\bf}
For $\alpha^*(t)= \min\alpha_D(t)$ and corresponding constants
$C_1$ and $C_2$ the following upper bound for the rate of convergence holds:

\begin{equation}
\|\vx(t)\| \le C_1 C_2 e^{-\int_0^t\alpha^*(\tau)\,
d\tau}\|\vx(0)\|. \label{05}
\end{equation}
\end{theorem}

\noindent
Note that if the matrix $A(t)$ is essentially non-negative then
the method based on differential inequalities
yields the same results as the method based
on the logarithmic norm. Thus the result of the {\it Theorem}~3 can also be
obtained using differential inequalities.

For some processes $X(t)$ belonging to classes (i)--(iv)
the method based on differential inequalities leads to such upper bounds,
which are better than those obtained using the both previous methods.
Several such settings are illustrated below.
Consider a homogeneous Markov chain $X(t)$ belonging to
class (iii) with the constant arrival intensity $\lambda$
and constant bulk service intensity $b_S(t)=b$ and $b_{k}(t)=0$, $1 \le k \le S-1$.
In this case both the method based on logarithmic norm
and the method based on Lyapunov functions do not yield any
upper bound, whereas with the
differential inequalities we can obtain a meaningful result.
Indeed, the matrix  $B^*$, given by (\ref{class3-fin}),
and the matrix  $B^{**}$ takes the following form:
\begin{equation} { B^* =  \left(
\begin{array}{cccccc}
-\lambda  & 0 & 0 & \cdots & 0 & - b \\
\lambda  & -\lambda & 0 & \cdots & 0 & - b \\
\ddots & \ddots & \ddots & \ddots & \ddots  & \ddots\\
0 & 0 & 0 & \cdots & \lambda & -\left(\lambda+ b\right)
\end{array}
\right)}, \label{class3*0}
\end{equation}
%%%%%%%%%%%%%%%%%%%
\begin{equation}
B^{**} = {\small \left(
\begin{array}{ccccccc}
-\lambda  & 0 & 0 & 0 &\cdots & 0 & -b \frac{d_{1}}{d_{S}} \\
 \lambda \frac{d_{2}}{d_{1}} & -\lambda & 0 & 0 &\cdots & 0 & -b \frac{d_{2}}{d_{S}} \\
0 & \lambda \frac{d_{3}}{d_{2}} & -\lambda & 0 & \cdots & 0 & -b \frac{d_{3}}{d_{S}} \\
\ddots & \ddots & \ddots & \ddots & \ddots & \ddots  & \ddots \\
0 & 0 & 0 & 0 & \cdots & -\lambda & - b  \frac{d_{S-1}}{d_{S}} \\
0 & 0 & 0 & 0 & \cdots & 0 & -\lambda - b
\end{array}
\right)}.
\end{equation}

\smallskip

\noindent Assume that $\{d_i, 1 \le i \le S\}$ are given and put $z_k(t)=d_kx_k(t)$. Then
we have

$$
\sum^{S}_{i=1} {dz_{i}(t) \over dt}
=
-\lambda
\sum^{S-1}_{i=1}
\left (1-\frac{d_{i+1}}{d_{i}} \right ) z_{i}(t)
-$$ $$-\left(\lambda+b \sum^{S}_{i=1} \frac{d_{i}}{d_{S}} \right)  z_{S}(t).
$$

%$$\sum^{S}_{i=1} z'_{i}(t)= -\lambda (1-\frac{d_{2}}{d_{1}})  z_{1}(t)-\lambda (1-\frac{d_{3}}{d_{2}})  z_{2}(t)-\lambda (1-\frac{d_{4}}{d_{3}})  z_{3}(t)-\cdots-$$
%$$-\left(\lambda+b \left(1+\frac{d_{1}}{d_{S}}+\frac{d_{2}}{d_{S}}+\cdots+ \frac{d_{S-1}}{d_{S}}\right)\right)  z_{S}(t).$$

Since $x_i(t)$ can be of different signs we have to consider all the
possible sign changes. It is convenient to start with the case when
there are no changes of signs. Let all $x_i(t)$ be positive. Put
$d_i=\varepsilon^i$, $1 \le i \le S$, for some $0<\varepsilon<1$.
Then
%%%%%%%%%%%%
%$$\sum^{S}_{i=1} z'_{i}(t)= -\lambda (1-\frac{d_{2}}{d_{1}})  z_{1}(t) -\lambda (1-\frac{d_{3}}{d_{2}})  z_{2}(t)-\lambda (1-\frac{d_{4}}{d_{3}})  z_{3}(t)-\cdots-$$
%$$
%-\left(\lambda+b \left(1+\frac{d_{1}}{d_{S}}+\frac{d_{2}}{d_{S}}+\cdots+
%\frac{d_{S-1}}{d_{S}}\right)\right)  z_{S}(t)=$$
%$$=-\lambda (1-\varepsilon)  z_{1}(t)-\lambda (1-\varepsilon)  z_{2}(t)-\lambda (1-\varepsilon)  z_{3}(t)-\cdots-$$
%$$
%-\left(\lambda +
%b \left(1+\frac{1}{\varepsilon}+\frac{1}{\varepsilon^2}+\cdots+
%\frac{1}{\varepsilon^{S-1}}\right)\right)  z_{S}(t),$$
%%%%%%%%%%%

$$
\sum^{S}_{i=1} {dz_{i}(t) \over dt}
=
-\lambda (1-\varepsilon  )
\sum^{S-1}_{i=1} z_{i}(t)
-\left(\lambda+b \sum^{S}_{i=1} \frac{1}{\varepsilon^{i-1}} \right)  z_{S}(t),
$$

\noindent and we obtain $\alpha=\lambda \left(1-\varepsilon\right)$.

The next case is when there is a single change of signs. Let all
$x_1(t),\dots,x_k(t)$ be positive for some $k$, $1 \le k \le S-1$,
and all $x_{k+1}(t),\dots,x_S(t)$ be negative. Put
$d_i=\varepsilon^{S-k+i}$, $1 \le i \le k$, and $d_{i}= -
\varepsilon^{i-k}$, $k+1 \le i \le S$. Then
%\begin{equation}
$$\sum^{S}_{i=1} z'_{i}(t)= -\lambda (1-\frac{d_{2}}{d_{1}})  z_{1}(t) -\lambda (1-\frac{d_{3}}{d_{2}})  z_{2}(t)-$$ $$-\lambda (1-\frac{d_{4}}{d_{3}})  z_{3}(t)-\cdots-$$
$$
-\left(\lambda+b \left(1+\frac{d_{1}}{d_{S}}+\frac{d_{2}}{d_{S}}+\cdots+
\frac{d_{S-1}}{d_{S}}\right)\right)  z_{S}(t)=$$
$$=-\lambda (1-\varepsilon)  z_{1}(t) -\lambda (1-\varepsilon)  z_{2}(t)-\lambda (1-\varepsilon)  z_{3}(t)-\cdots-$$
$$
-\lambda (1-\varepsilon)
z_{k-1}(t) -\lambda (1+\frac{1}{\varepsilon^{S-1}})
z_{k}(t) -\cdots-
$$
$$
-\left(\lambda + b
\left(1-\varepsilon-\varepsilon^{2}-...-\varepsilon^{k}+\frac{1}{\varepsilon^{S-k-1}}+\right.\right.$$ $$\left.\left.+\frac{1}{\varepsilon^{S-k-2}}+\cdots+
\frac{1}{\varepsilon}\right)\right)  z_{S}(t) \le$$
$$
\le -\lambda\left(1-\varepsilon\right)\sum^{S}_{i=1} z_{i}(t),
$$
 and we have that $\alpha=\lambda \left(1-\varepsilon\right)$.

Now consider the case when there are exactly two changes of signs.
Let all $x_1(t),\dots, x_k(t)$ be positive for some $1 \le k \le
S-2$, all $x_{k+1}(t),\dots,x_s(t)$ be negative for some $k+1 \le s
\le S-1$ and all $x_{s+1}(t),\dots, x_S(t)$ be positive. Let
$d_i=\varepsilon^{S-k+i}$ for $1 \le i \le k$, $d_{i} = -
\varepsilon^{S-s-k+i}$ for $k+1 \le i \le s$ and $d_{i} =
\varepsilon^{i-s}$, for $s+1 \le i \le S$. We have
%\begin{equation}
{\small
$$\sum^{S}_{i=1} z'_{i}(t)= -\lambda (1-\frac{d_{2}}{d_{1}})  z_{1}(t)-$$ $$-\lambda (1-\frac{d_{3}}{d_{2}})  z_{2}(t)-\lambda (1-\frac{d_{4}}{d_{3}}) z_{3}(t)-\cdots-$$
$$
-\left(\lambda+b \left(1+\frac{d_{1}}{d_{S}}+\frac{d_{2}}{d_{S}}+\cdots+
\frac{d_{S-1}}{d_{S}}\right)\right)  z_{S}(t)=$$
$$=-\lambda (1-\varepsilon)  z_{1}(t)-\lambda (1-\varepsilon) z_{2}(t)-\lambda (1-\varepsilon)  z_{3}(t)-\cdots-$$
$$
-\lambda (1-\varepsilon)
z_{k-1}(t)-\lambda (1+\frac{1}{\varepsilon^{s-1}})
z_{k}(t)-\lambda (1-\varepsilon)
z_{k+1}(t)-\cdots-
$$
$$
-\lambda (1-\varepsilon)  z_{s-1}-\lambda (1+\frac{1}{\varepsilon^{S-k}})
z_{s}(t)-\lambda (1-\varepsilon) z_{s+1}(t)-\cdots-
$$
$$ -\left(\lambda +
b \left(1+\varepsilon^{s-k+1}+\varepsilon^{s-k+2}+...+\varepsilon^{s}-\right.\right. $$
$$\left.\left.-\varepsilon-\varepsilon^2-...-\varepsilon^{s-k}+
\frac{1}{\varepsilon^{S-k-1}}+\right.\right. $$
$$\left.\left. +\frac{1}{\varepsilon^{S-k-2}}+\cdots+
\frac{1}{\varepsilon}\right)\right)  z_{S}(t)\le -\lambda\left(1-\varepsilon\right)\sum^{S}_{i=1} z_{i}(t),$$}
%\end{equation}

\noindent and $\alpha=\lambda \left(1-\varepsilon\right)$. Note that
the total number of sign changes does not exceed $S-1$. On each
change of sign when going from $x_s(t)$ to $x_{s+1}(t)$ we put
$d_{s+1}$ equal to $\varepsilon^{S-m+1}$, where $m$ is the number of
the last element in the current period of consistency (i.e. when
there is no change of signs). Then eventually we arrive at the
following upper bound $\|\vx(t)\| \le C_1   C_2  e^{-\lambda
\left(1-\varepsilon\right) t}\|\vx(0)\|$, with $C_1 C_2 =
\varepsilon^{1- S}$.

The following example shows how
the method based on differential inequalities
can be applied for inhomogeneous Markov chains with finite state space.

\bigskip

{\bf Example 2.} Consider the Markov process $X(t)$ that describes
the evolution of the total number of customers in the $M(t)/M^X(t)/1/S$ queue with bulk services,
when all transition intensities are periodic functions of time.
Let the arrival intensities
be $\lambda_k(t)
=\lambda(t)=10\left(2+ \sin (2 \pi t)\right)$,
and the service intensities be
$b_k(t)=0$, for $1 \le k < S$, and $b_S(t)=m^{-2}\left(2+\cos 2\pi t\right)$
for some $ m \ge 1$.
Such $X(t)$ belongs to class (iii).
The matrix $B^{**}$ for such $X(t)$ has the form $B^{**}(t)=(b_{ij}^{**}(t))$, where
%$$b_{i,i+1}^{**}(t)=0,$$
$$b_{i,i}^{**}(t)=-10(2+\sin(2\pi t)),$$
$$b_{i+1,i}^{**}(t)=10(2+\sin(2\pi t)) \frac{d_{i+1}}{d_{i}},$$
$$b_{i,S}^{**}(t)=-m^{-2}\left(2+\cos(2\pi t)\right) \frac{d_{i}}{d_{S}}, ~~  i< S,$$
$$b_{SS}^{**}(t)=-10(2+\sin(2\pi t)) - m^{-2}\left(2+\cos(2\pi t)\right), ~~ i=S.$$

%\begin{equation*}
%b_{i,i}^{**}(t)=
%\left\{ \!
%\begin{array}{cc}
%b_{i,S}^{**}(t)=-m^{-2}\left(2+\cos(2\pi t)\right) \frac{d_{i}}{d_{S}}, & \mbox {if } i\neq S,
%\\ b_{SS}^{**}(t)=-10(2+\sin(2\pi t)) - m^{-2}\left(2+\cos(2\pi t)\right), & \mbox {if } i=S,
%\end{array}
%\right.
%\end{equation*}

Then for any initial condition ${\bf \vx}(0)$
we can deduce the following two upper bounds on the rate of convergence:
$$
\|\vx(t)\| \le \varepsilon^{1- S}   e^{-{ \int^{t}_{0}\left(1-\varepsilon\right) \lambda(\tau) d\tau}}\|\vx(0)\|,
$$
$$
\
\|\vx(t)\| \le \varepsilon^{1- S}   e^{-10\left(1-\varepsilon\right)  t}\|\vx(0)\|.
$$

These bounds are not sharp (the leftmost is better among the two)
but the other two methods give essentially worse results. As in the
{\it Example~1}, these bounds can be used in the approximation of
the limiting distribution of $X(t)$. For example, let $S=40$ and
$m=1$. In Figs.~9--16 we can see the behaviour of the conditional
expected number $E(X(t)|X(0))$ of customers in the queue at instant
$t$ and the state probabilities $p_0(t)$, $p_{20}(t)$ and
$p_{40}(t)$ as functions of time $t$ under different initial
conditions $X(0)$.

We conclude the section by emphasizing that the method of differential inequalities
may lead to meaningful upper bounds for the rate of convergence
even in the case of countable state space $\mathcal{X}$.
For example, consider a homogeneous countable (i.e. $S=\infty$) Markov process
$X(t)$ belonging to class (iii) with constant arrival
intensities $\lambda$
and batch service intensities $b_2(t)=\mu >0$ and $b_k(t)=0$
for $k\neq 2$. Hence the matrix $B^*$, given by (\ref{class3-fin}),
takes the form
\begin{equation}
B^* = {\scriptsize  \left(
\begin{array}{cccccc}
-\lambda  &  - \mu
 & \mu & \cdots & \cdots \\
\lambda & -\left(\lambda + \mu\right) & 0 & \mu & \cdots\\
\vdots & \ddots & \ddots & \ddots & \ddots  \\
\vdots & \cdots & \cdots & \lambda &
-\left(\lambda+ \mu\right) \cdots \\
\cdots & \cdots & \cdots & \cdots & \cdots & \cdots \\
\end{array}
\right)}.\label{class3-inf-ineq}
\end{equation}
 \noindent
In such case, to our best knowledge, the method of differential
inequalities is the only method, with which we can obtain the
ergodicity of the chain and explicit estimates of the rate of
convergence (see the details in \cite{Satin2019}).

%%%%%%% case (i)

%\begin{center}
%\begin{figure}[!ht]
%\begin{center}
%  \includegraphics[scale=0.35]{4-l20phi.eps}
%\end{center}
%\caption{Case (i), the arrival intensity is $\lambda^*(t;20)$. Approximation of the limiting
%mean number $\varphi (t)$ of customers in the system for $t\in[5,6]$.}
%\label{fig:43}       % Give a unique label
%\end{figure}
%\end{center}

\begin{center}
\begin{figure}[!ht]
\begin{center}
  \includegraphics[scale=0.35]{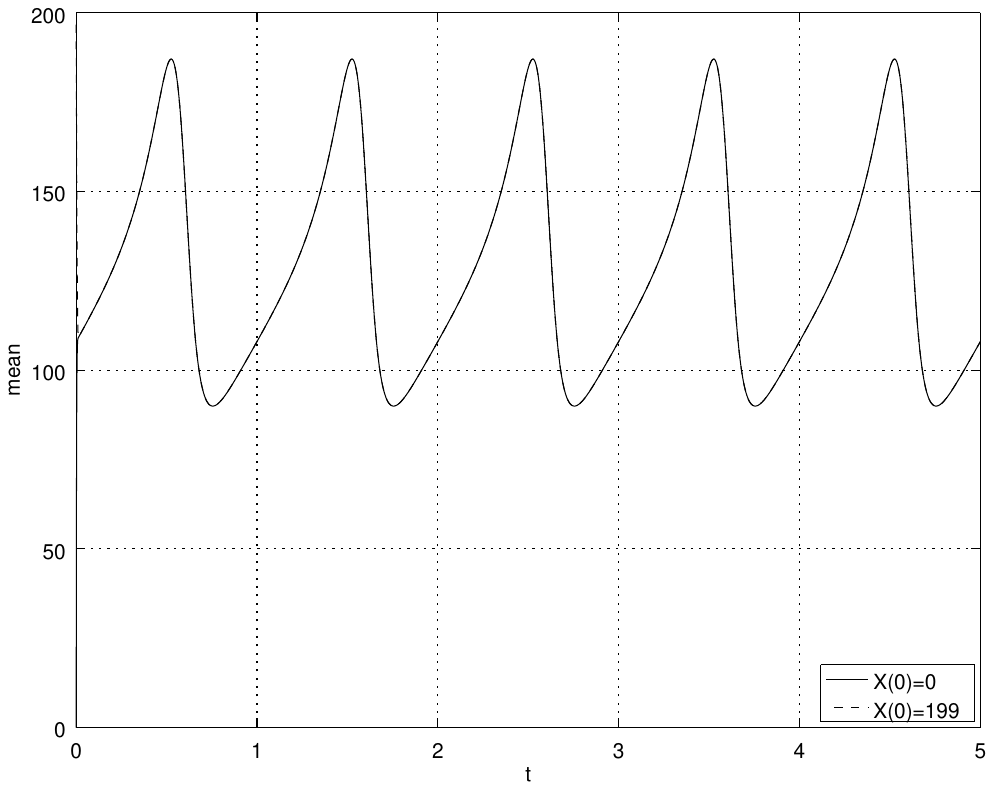}%\vspace{-3cm}
\end{center}
\caption{ Example 1.  The expected number $E(X(t)|X(0))$ of customers in the queue for $t\in[0,5]$ with the initial
condition $X(0)=0$.}
\label{fig:01}       % Give a unique label
\end{figure}
\end{center}

\begin{center}
\begin{figure}[!ht]
\begin{center}
  \includegraphics[scale=0.35]{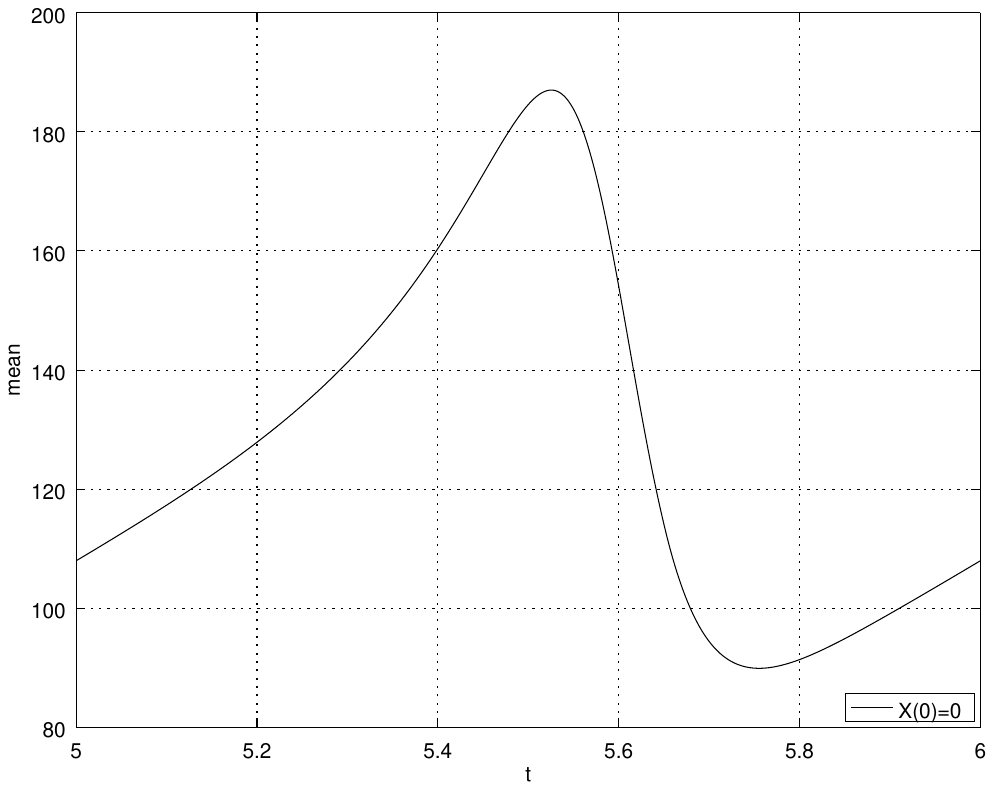}%\vspace{-3cm}
\end{center}
\caption{ Example 1.  The expected number $E(X(t)|X(0))$ of customers in the queue for $t\in[5,6]$ with the initial
condition $X(0)=0$.}
\label{fig:03}       % Give a unique label
\end{figure}
\end{center}

\begin{center}
\begin{figure}[!ht]
\begin{center}
  \includegraphics[scale=0.35]{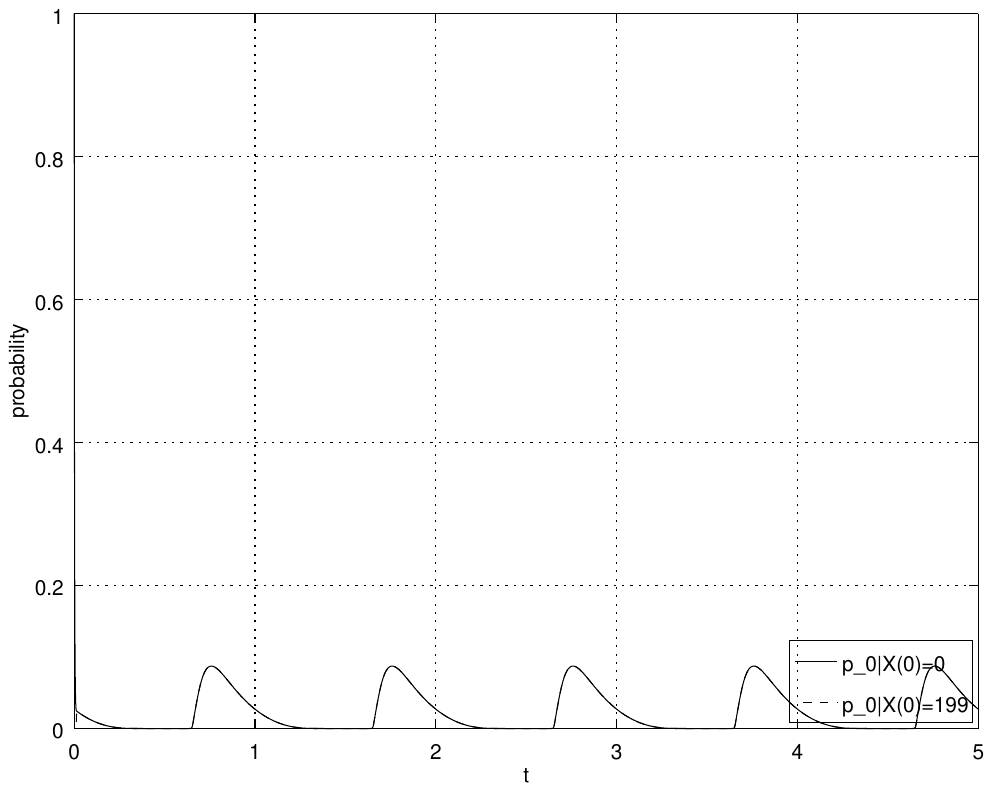}%\vspace{-3cm}
\end{center}
\caption{ Example 1. The probability $p_0(t)$ of empty queue for
$t\in[0,5]$ with the initial condition $X(0)=0$.}
\label{fig:04}       % Give a unique label
\end{figure}
\end{center}

\begin{center}
\begin{figure}[!ht]
\begin{center}
  \includegraphics[scale=0.35]{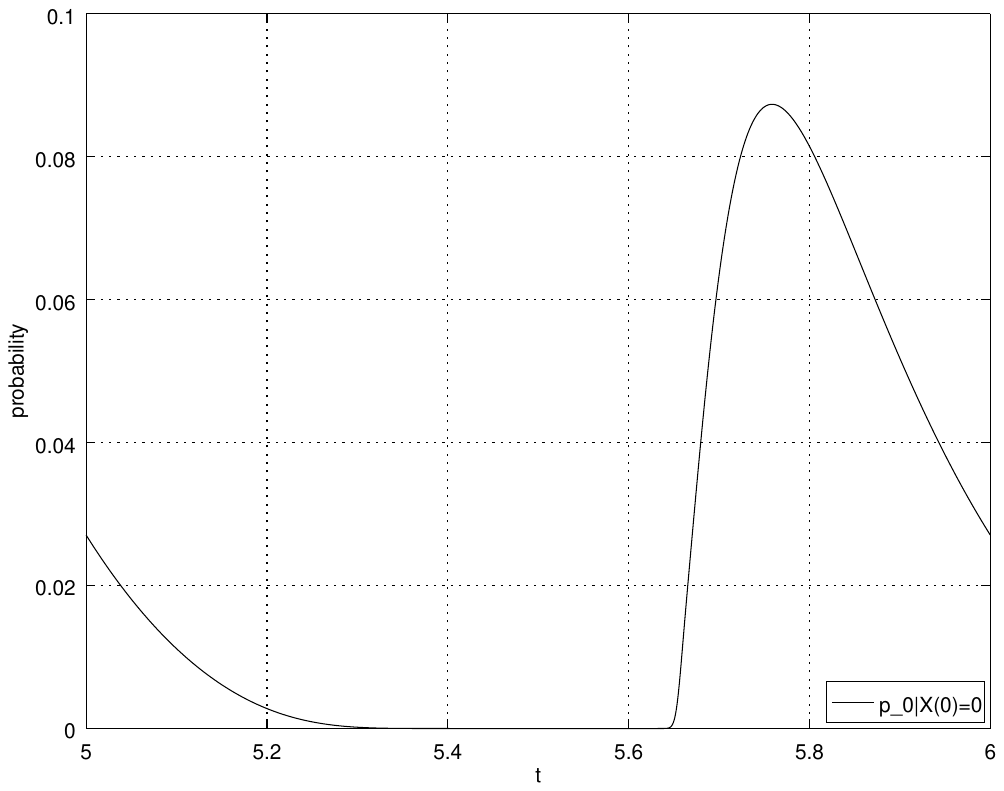}%\vspace{-3cm}
\end{center}
\caption{ Example 1. The probability $p_0(t)$ of empty queue  for
$t\in[5,6]$ with the initial condition $X(0)=0$.}
\label{fig:05}       % Give a unique label
\end{figure}
\end{center}

\begin{center}
\begin{figure}[!ht]
\begin{center}
  \includegraphics[scale=0.35]{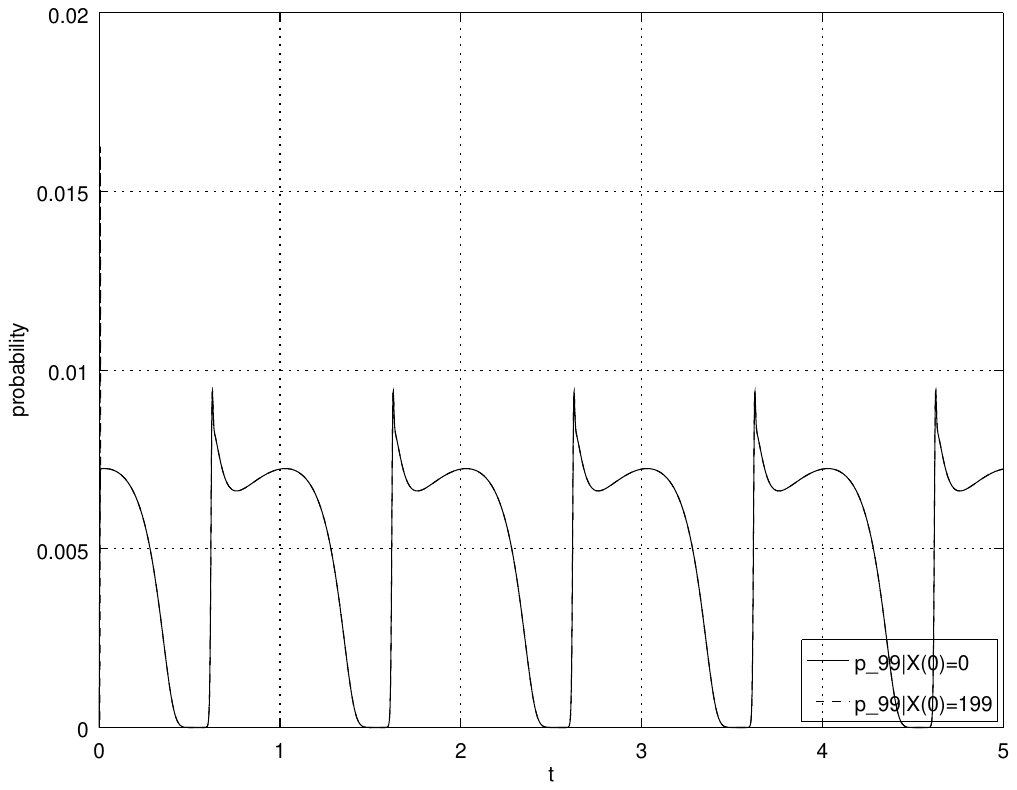}%\vspace{-3cm}
\end{center}
\caption{ Example 1. The probability $p_{99}(t)$ for $t\in[0,5]$ with the
initial condition $X(0)=0$.}
\label{fig:06}       % Give a unique label
\end{figure}
\end{center}

\begin{center}
\begin{figure}[!ht]
\begin{center}
  \includegraphics[scale=0.35]{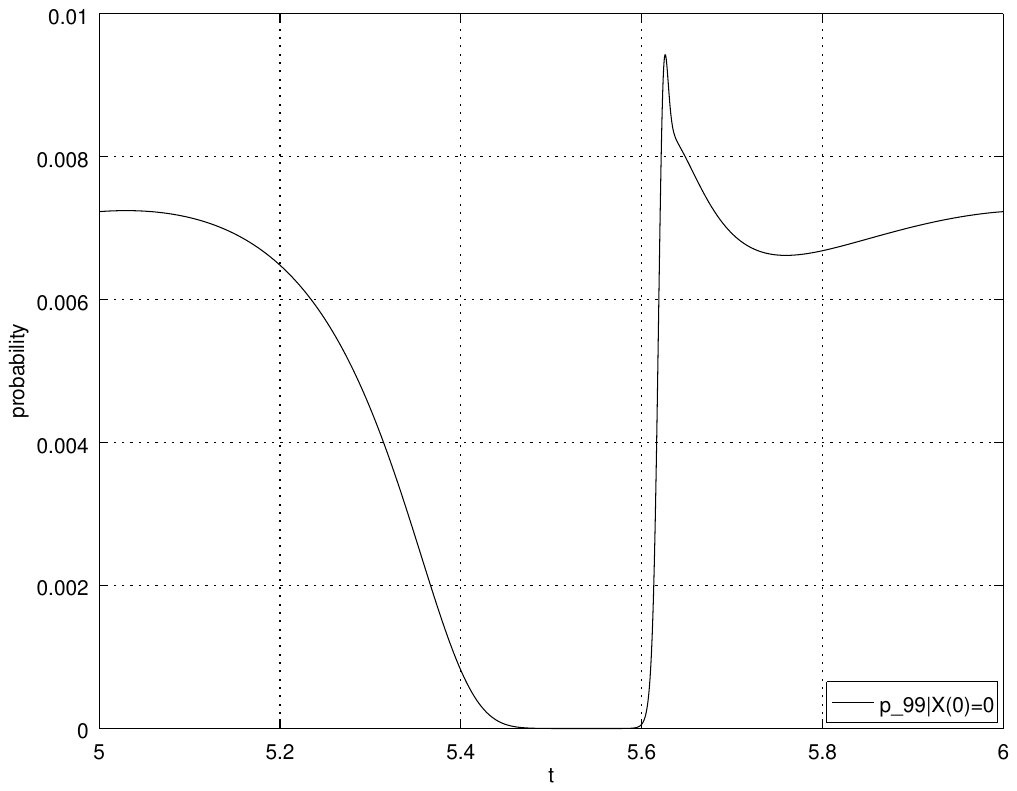}%\vspace{-3cm}
\end{center}
\caption{ Example 1. The probability  $p_{99}(t)$ for $t\in[5,6]$ with the
initial condition $X(0)=0$.}
\label{fig:07}       % Give a unique label
\end{figure}
\end{center}

\begin{center}
\begin{figure}[!ht]
\begin{center}
  \includegraphics[scale=0.35]{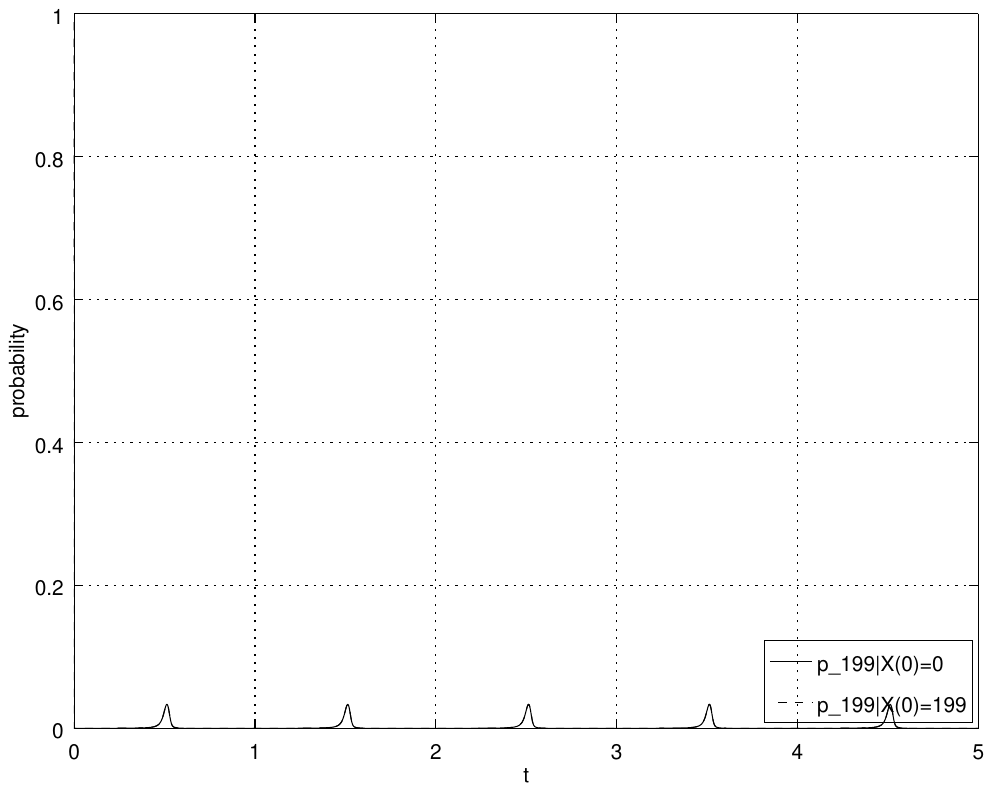}%\vspace{-3cm}
\end{center}
\caption{ Example 1. The probability  $p_{199}(t)$ for $t\in[0,5]$ with the
initial condition $X(0)=0$.}
\label{fig:08}       % Give a unique label
\end{figure}
\end{center}

\begin{center}
\begin{figure}[!ht]
\begin{center}
  \includegraphics[scale=0.35]{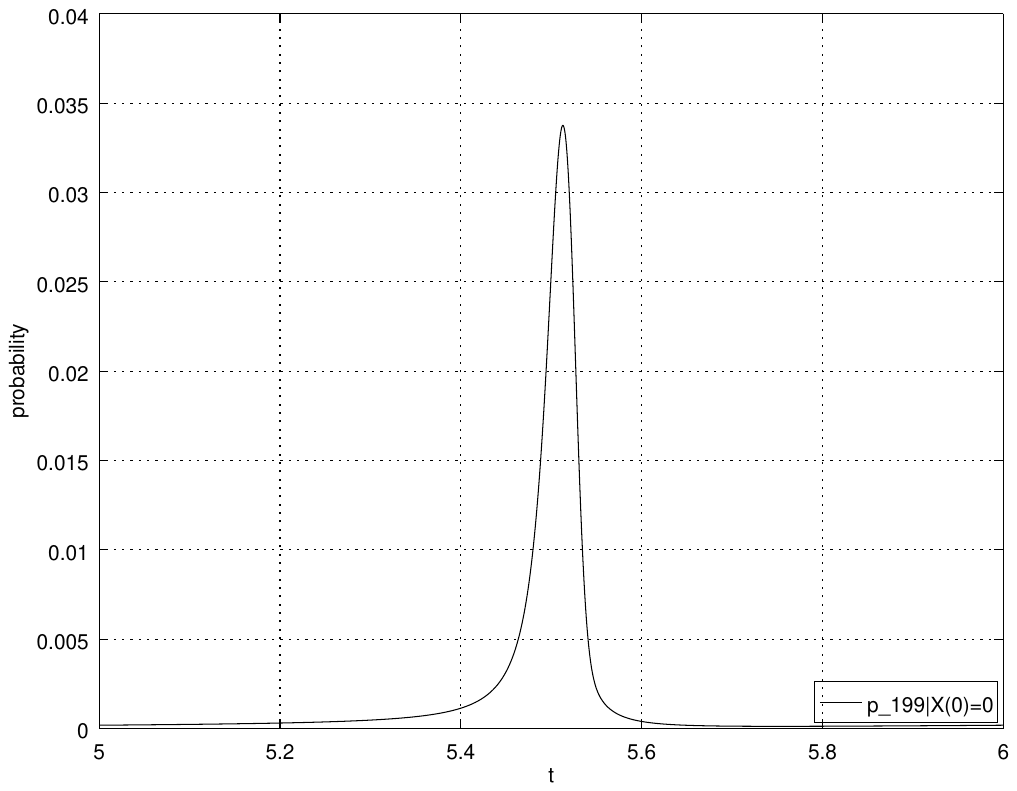}%\vspace{-3cm}
\end{center}
\caption{ Example 1. The probability  $p_{199}(t)$ for $t\in[5,6]$ with the
initial condition $X(0)=0$.}
\label{fig:09}       % Give a unique label
\end{figure}
\end{center}

\begin{center}
\begin{figure}[!ht]
\begin{center}
  \includegraphics[scale=0.35]{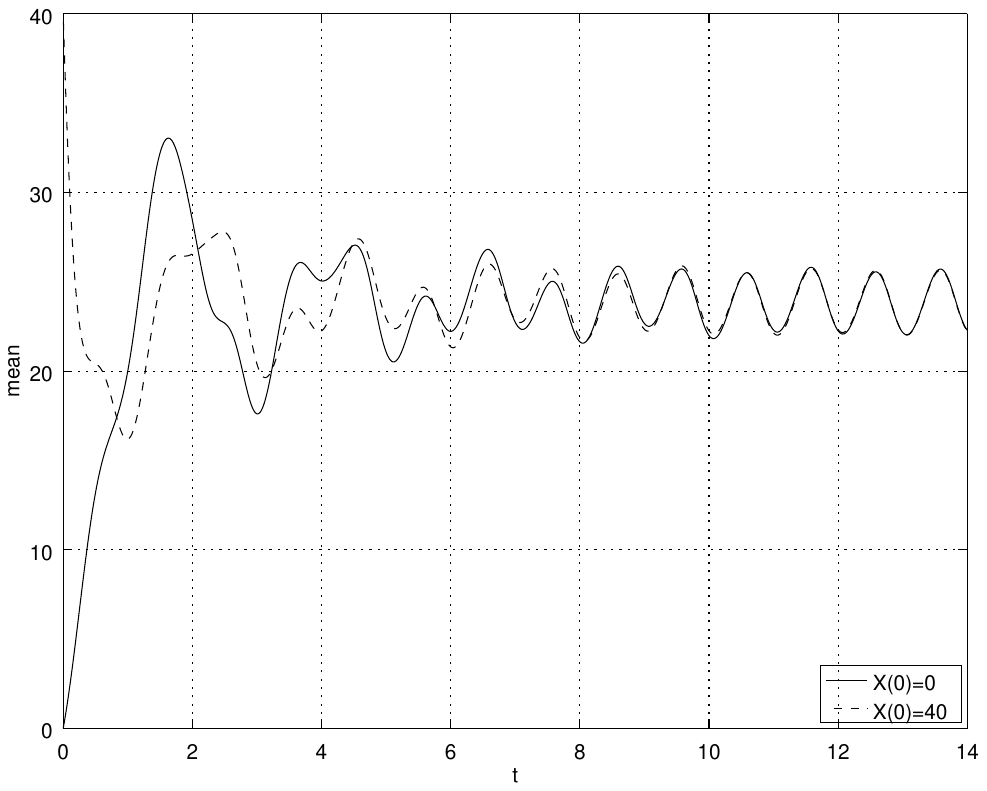}%\vspace{-3cm}
\end{center}
\caption{ Example 2.  The expected number $E(X(t)|X(0))$ of customers in the queue for $t\in[0,14]$ with the initial conditions $X(0)=0$ and $X(0)=S$.}
\label{fig:10}       % Give a unique label
\end{figure}
\end{center}

\begin{center}
\begin{figure}[!ht]
\begin{center}
  \includegraphics[scale=0.35]{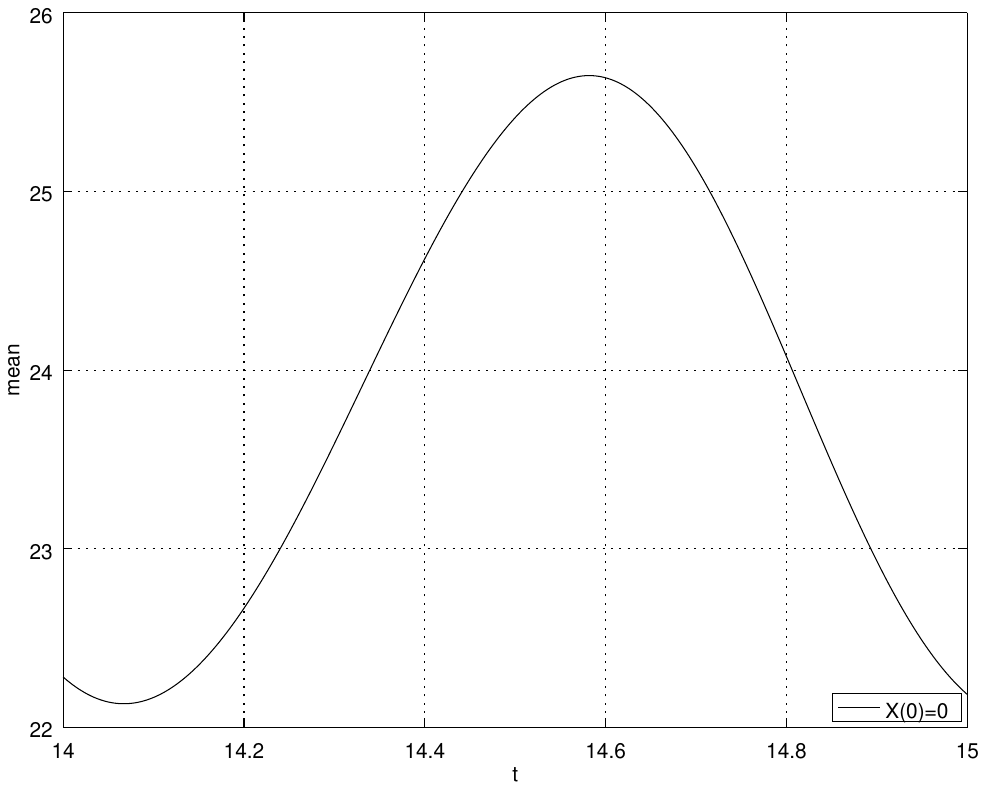}%\vspace{-3cm}
\end{center}
\caption{ Example 2.  The expected number $E(X(t)|X(0))$ of customers in the queue for $t\in[14,15]$ with
initial condition $X(0)=0$.}
\label{fig:11}       % Give a unique label
\end{figure}
\end{center}

\begin{center}
\begin{figure}[!ht]
\begin{center}
  \includegraphics[scale=0.35]{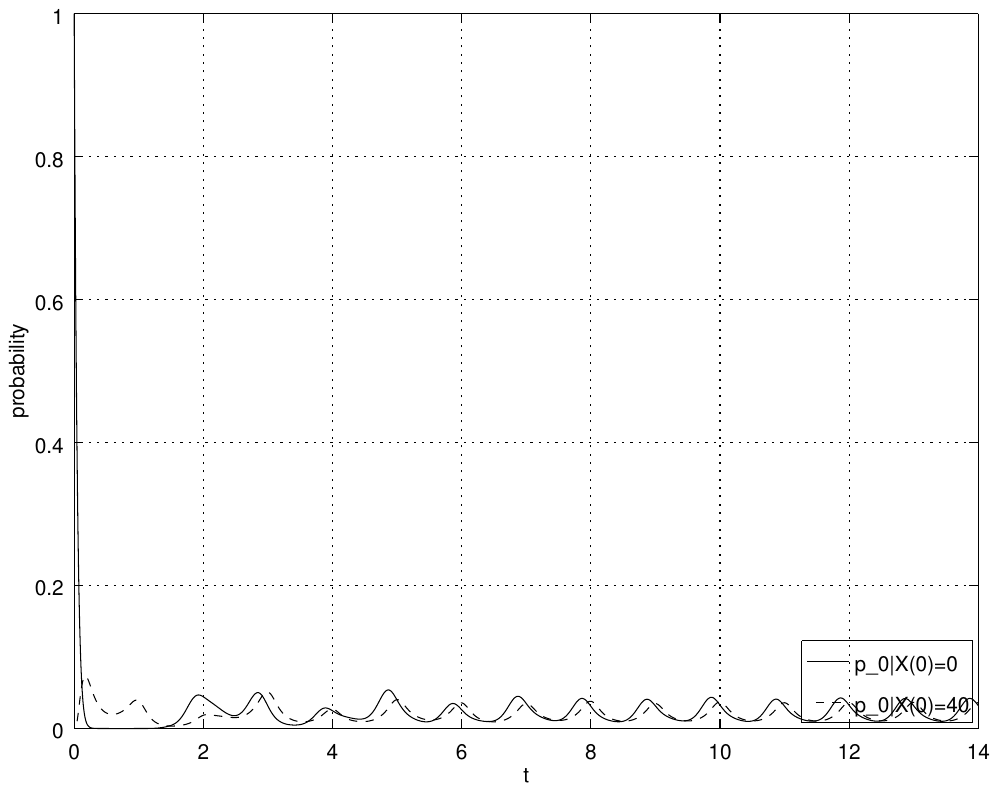}%\vspace{-3cm}
\end{center}
\caption{ Example 2. The probability $p_0(t)$ of empty queue for
$t\in[0,14]$ with the initial conditions $X(0)=0$ and $X(0)=S$.}
\label{fig:12}       % Give a unique label
\end{figure}
\end{center}

\begin{center}
\begin{figure}[!ht]
\begin{center}
  \includegraphics[scale=0.35]{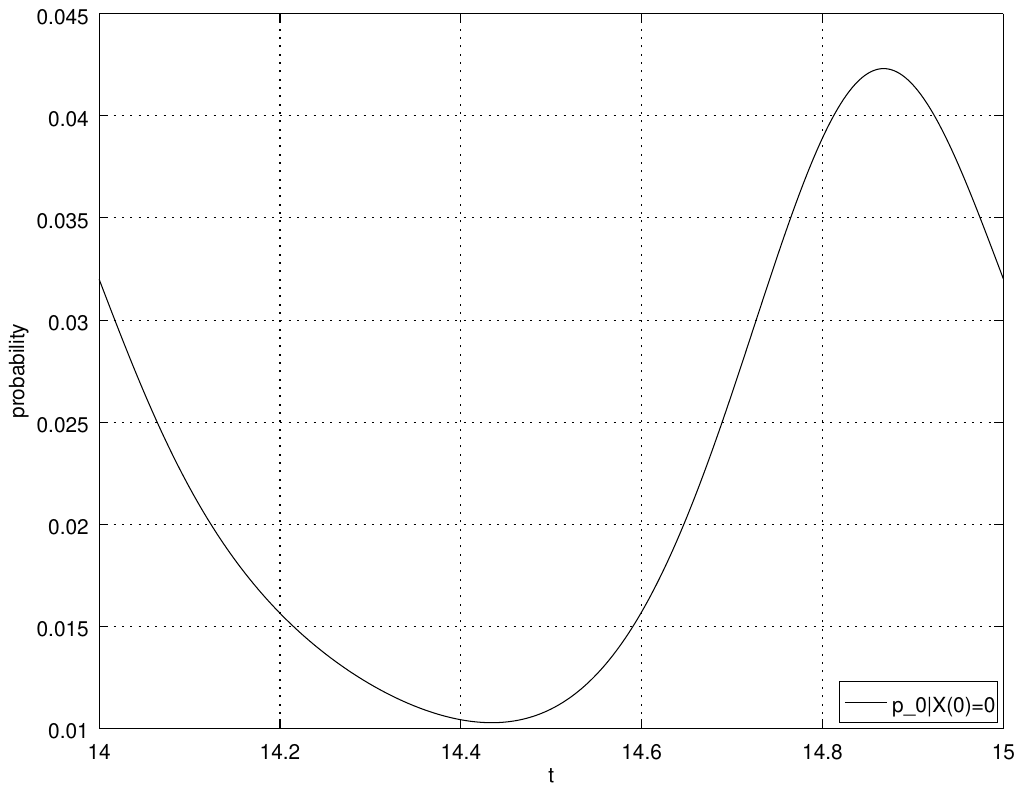}%\vspace{-3cm}
\end{center}
\caption{ Example 2. The probability $p_0(t)$ of empty queue  for
$t\in[14,15]$ with the initial condition $X(0)=0$.}
\label{fig:13}       % Give a unique label
\end{figure}
\end{center}

\begin{center}
\begin{figure}[!ht]
\begin{center}
  \includegraphics[scale=0.35]{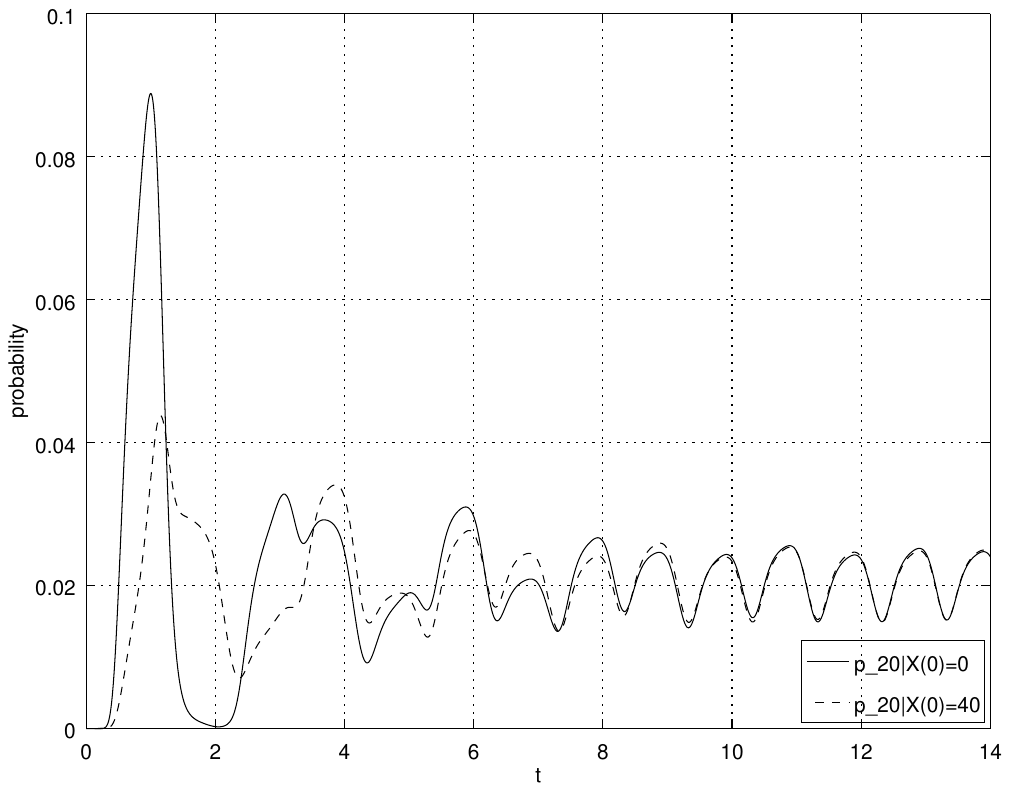}%\vspace{-3cm}
\end{center}
\caption{ Example 2. The probability  $p_{20}(t)$ for $t\in[0,14]$ with the
initial conditions $X(0)=0$ and $X(0)=S$.}
\label{fig:14}       % Give a unique label
\end{figure}
\end{center}

\begin{center}
\begin{figure}[!ht]
\begin{center}
  \includegraphics[scale=0.35]{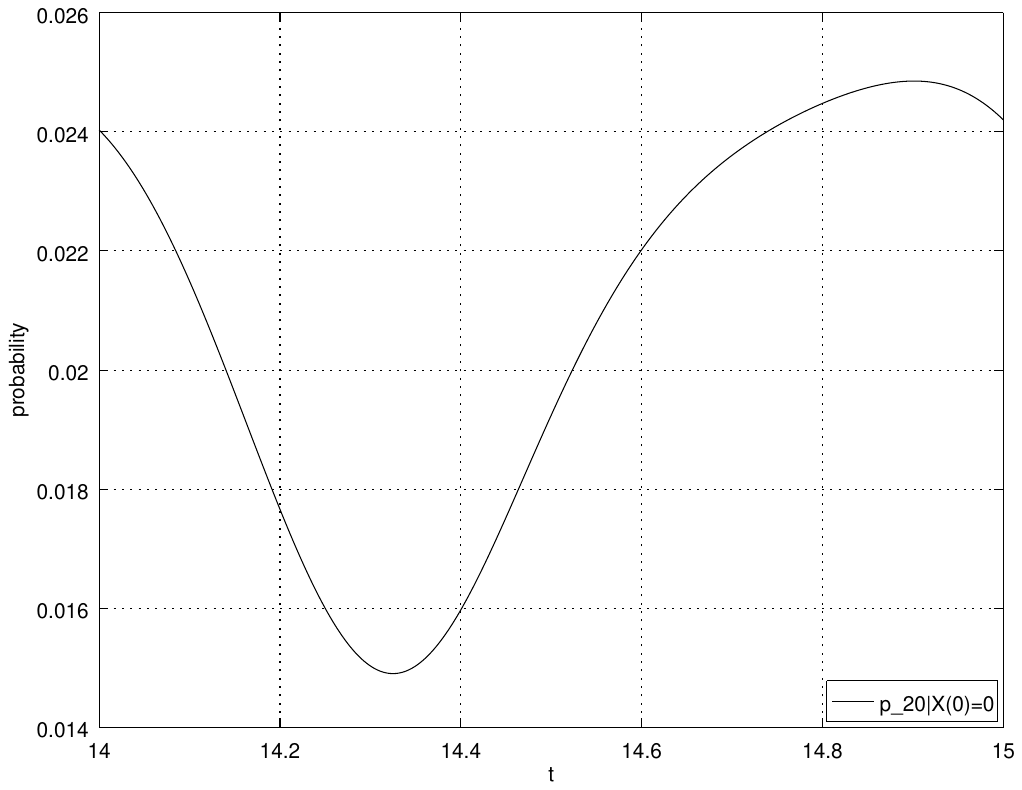}%\vspace{-3cm}
\end{center}
\caption{ Example 2. The probability $p_{20}(t)$ for $t\in[14,15]$ with the
initial condition $X(0)=0$.}
\label{fig:15}       % Give a unique label
\end{figure}
\end{center}

\begin{center}
\begin{figure}[!ht]
\begin{center}
  \includegraphics[scale=0.35]{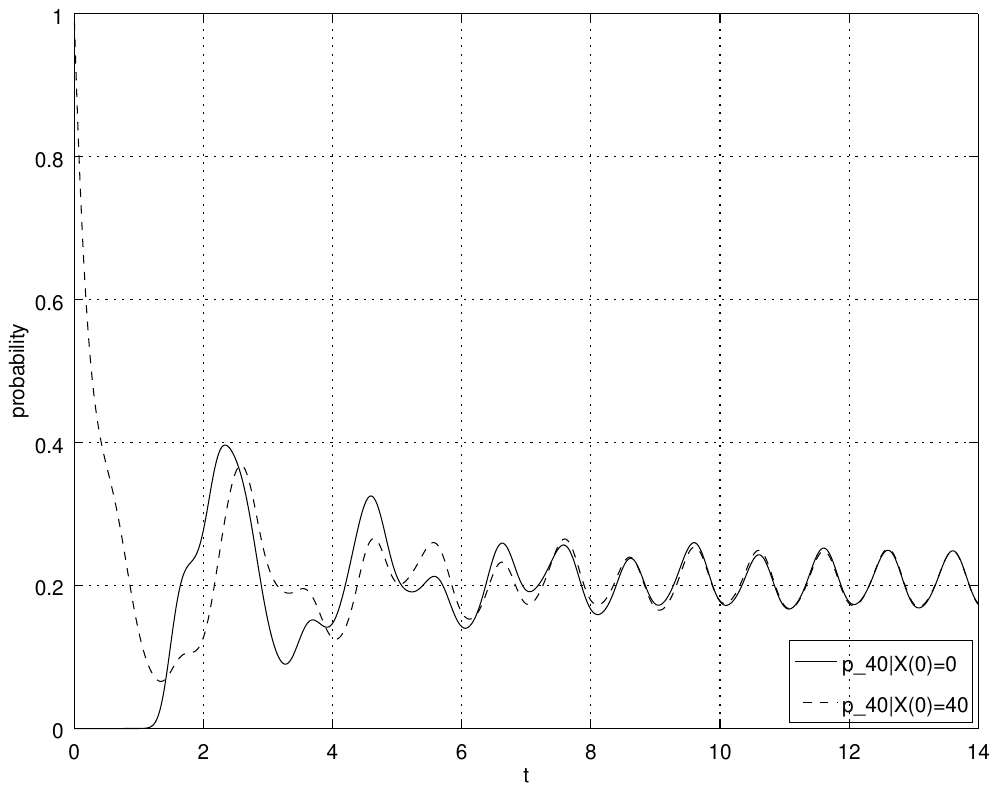}%\vspace{-3cm}
\end{center}
\caption{ Example 2. The probability  $p_{40}(t)$ for $t\in[0,14]$ with the
initial conditions $X(0)=0$ and $X(0)=S$.}
\label{fig:16}       % Give a unique label
\end{figure}
\end{center}

\begin{center}
\begin{figure}[!ht]
\begin{center}
  \includegraphics[scale=0.35]{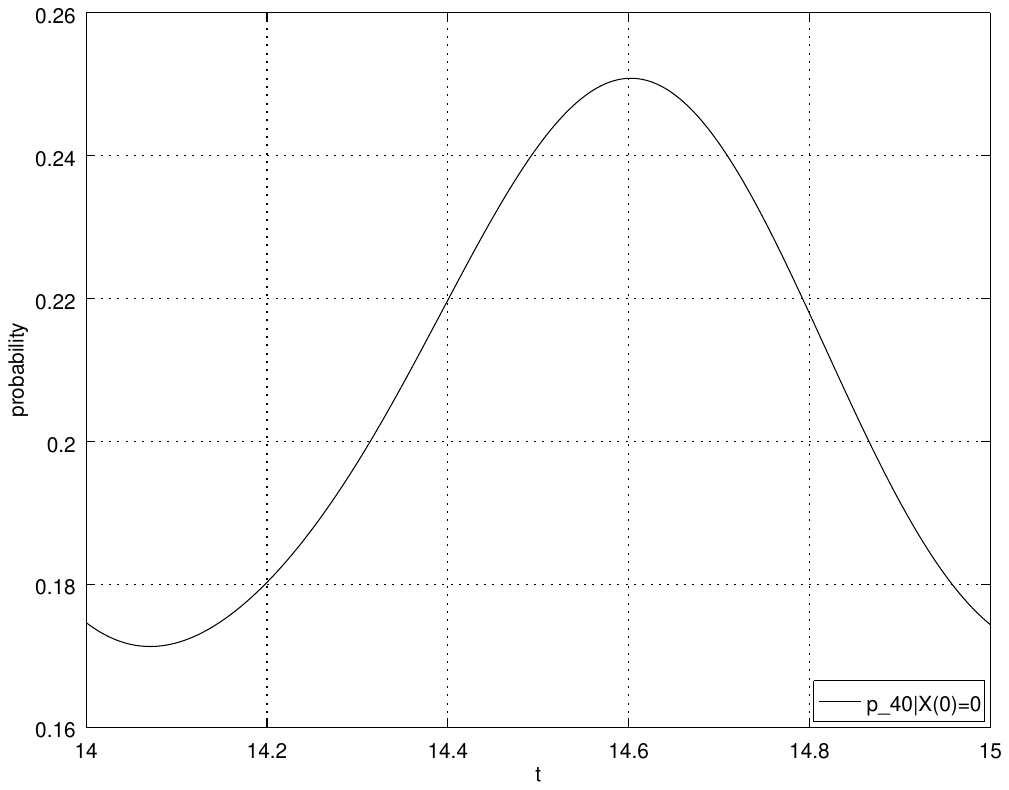}%\vspace{-3cm}
\end{center}
\caption{ Example 2. The probability  $p_{40}(t)$ for $t\in[14,15]$ with the
initial condition $X(0)=0$.}
\label{fig:09}       % Give a unique label
\end{figure}
\end{center}

\clearpage

\section{Conclusion}

The three methods considered in this paper provide various alternatives for the computation of the upper bounds for the rate of convergence to the limiting regime of (in)homogeneous continuous-time Markov processes.
Yet even for the four considered classes (i)--(iv) of Markov processes a single unified framework cannot be suggested: special cases do exist when none of the methods works well.

\section{Acknowledgment}
This research was supported by Russian Science Foundation under grant 19-11-00020.

\end{document}